\documentclass[11pt,reqno]{amsart}
\usepackage{graphicx}
\usepackage{verbatim}
\usepackage{textcomp}
\usepackage{amssymb}
\usepackage{cite}
\usepackage{amsmath}
\usepackage{latexsym}
\usepackage{amscd}
\usepackage{amsthm}
\usepackage{mathrsfs}
\usepackage{xypic}
\usepackage{bm}
\usepackage{url}
\usepackage{hyperref}

\newtheorem{thm}{Theorem}[section]

\newtheorem{lem}[thm]{Lemma}

\theoremstyle{definition}
\newtheorem{defn}{Definition}[section]

\theoremstyle{remark}
\newtheorem{rem}{Remark}[section]
\numberwithin{equation}{section} 

\begin{document}
\title{Isoparametric hypersurfaces in Minkowski spaces}
\author{Qun He, Songting Yin}
\address{Department of Mathematics, Tongji University, Shanghai, 200092, China}
 \email{hequn$@$tongji.edu.cn}

\address{Department of Mathematics and Computer Science, Tongling University,
Tongling, 244000, China} \email{yst419$@$163.com}
\author{Yibing Shen}
\address{Department of Mathematics, zhejiang University, Hangzhou, 310028 zhejiang, China}
 \email{yibingshen$@$zju.edu.cn}
\thanks{Project supported by NNSFC(no.11471246) and  NSFHE(no.KJ2014A257)}
\maketitle

\begin{abstract}
In this paper, we introduce isoparametric functions and
isoparametric hypersurfaces in Finsler manifolds and give the
necessary and sufficient conditions for a transnormal function to be
isoparametric.  We then prove that hyperplanes, Minkowski
hyperspheres and $F^*$-Minkowski cylinders in a Minkowski space with
$BH$-volume (resp. $HT$-volume) form are all isoparametric
hypersurfaces with one and two distinct constant principal
curvatures respectively. Moreover, we give a complete classification
of isoparametric hypersurfaces in Randers-Minkowski spaces and
construct a counter example, which shows that Wang's Theorem B in
\cite{WQ} does not hold in Finsler geometry.
\end{abstract}

\section{\textbf{Introduction}}

In Riemannian geometry, isoparametric hypersurfaces are a class of
important submanifolds studied by many geometers. Let $(M,g)$ be a
connected complete Riemannian manifold. An \emph{isoparametric
hypersurface} in $(M,g)$ is a regular level hypersurface of an {\it
isoparametric function} $f$, which satisfies
\begin{equation}\label{1.0}
\left\{\begin{aligned}
&|\nabla f|^2=\tilde{a}(f),\\
&\Delta f=\tilde{b}(f),
\end{aligned}\right.\end{equation}
where $\Delta$ denotes the Laplacian on $(M,g)$, $\tilde a$ is a
smooth function and $\tilde b$ is a continuous function on $M$. A
smooth function $f$ satisfying only the first equation of
(\ref{1.0}) is called a {\it transnormal function}. Early in 1938,
B. Segre (\cite{27}) proved that an isoparametric hypersurface in
$\mathbb{R}^n$ is either a hyperplane, a hypersphere, or a cylinder
$\mathbb{S}^{m-1}\times \mathbb{R}^{n-m}$. Later on, E. Cartan
(\cite{Ca}) showed that in a space form $M^n(c)$, a transnormal
function $f$ is isoparametric if and only if each regular level
hypersurface of $f$ has constant mean curvature, or equivalently, if
and only if each regular level hypersurface has constant principal
curvatures. If $k_1,k_2,\cdots,k_{g}$ are the all distinct principal
curvatures, then we have the following Cartan formula
\begin{align} \sum_{i\neq j}m_i\frac{c+k_jk_i}{k_j-k_i}=0,~~~~j=1,\ldots,g,
\label{1.1}\end{align} where $m_i$ is the  multiplicity of $k_i$.
Subsequent research focused  mainly on isoparametric hypersurfaces
in the $n$-sphere $\mathbb{S}^{n}$ and many important results have
been obtained(\cite{CR}\cite{QT}\cite{WQ}). Recently, Z. Tang and
his students gave a lot of fresh results on isoparametric
hypersurfaces in $\mathbb{S}^{n}$ (\cite{TY}\cite{QT}).

In Finsler  geometry, however, there are no any studies on
Finslerian isoparametric hypersurfaces so far. The major
difficulties are that there is no unified volume measure and  there
are many ways to define the induced volume form and the mean
curvature of a hypersurface in a Finsler manifold $(M,F)$. Usually,
the Busemann-Hausdorff volume  ($BH$-volume for short) form and the
Holmes-Thompson volume ($HT$-volume for short) form are used
(\cite{SZ}). Given a volume form $d\mu$ on $(M,F)$, we can define a
non-linear Finsler-Laplacian
$$\Delta f=\Delta^{\nabla f} f={\text{div}}({\nabla f})$$ with respect to $d\mu$, where $\nabla f$ denotes the gradient of
$f$ defined by means of the Legendre transformation. The  non-linear
Finsler-Laplacian $\Delta f$ is a well-known and very important
operator in Finsler geometry, on which there are many important
results(\cite{SZ}\cite{OS}\cite{YHS}, etc). We then can define the
isoparametric function and the isoparametric hypersurface in a
Finsler manifold $(M,F,d\mu)$ in the same manner as in Riemannian
geometry. In other words, a function  $f$ on $(M,F,d\mu)$ is said to
be \textit{isoparametric} if  there are two functions $\tilde a (t)$
and $\tilde b (t)$ such that $f$ satisfies (see Definition 4.1 for
details)
\begin{equation}\label{1.3} \left\{\begin{aligned}
&F(\nabla f)=\tilde{a}(f),\\
&\Delta f=\tilde{b}(f).
\end{aligned}\right.\end{equation}
A function $f$ satisfying the first equation of (1.3) is said to be
{\it transnormal}.

It appears that the equation (1.3) and (1.1) are similar in form.
But in the non Riemannian case, the gradient of a function and its
Laplacian are nonlinear. In general case, they cannot be expressed
explicitly. Even in the most special case, Randers-Minkowski space,
the expressions of the gradient and Laplacian are also quite
complicated, so that many methods of Riemannian geometry are no
longer adequate. One has to find some new methods and techniques for
handling equation (1.3).

The purpose of this paper is to study isoparametric functions and
isoparametric hypersurfaces in a Finsler manifold $(M,F,d\mu)$,
particularly, in a Minkowski space which is a generalization of the
Euclidean space. We shall prove the following results.

\begin{thm}\label{thm1-1} On an $n$-dimensional Finsler manifold $(M,F,d\mu)$, a transnormal function  $f$ is
isoparametric if and only if each regular level hypersurface  $N_t$
of  $f$ has constant $d\mu_{\textbf{n}}$-mean curvature
${H}_{\textbf{n}}$, where $\textbf{n}=\frac{\nabla  f }{F(\nabla  f
)}$. Particularly, if $M$ has constant flag curvature and constant
{\bf S}-curvature, then a transnormal function  $f$ is isoparametric
if and only if  all the principal curvatures of $N_t$ are constant.
\end{thm}

The above theorem can be viewed as a generalization of the
corresponding result in \cite{Ca}. As well known, a Minkowski space
is a vector space endowed with Minkowski metric which is Euclidean
metric without quadratic restriction. It is  natural to generalize
the results in Euclidean spaces to Minkowski spaces.
 We consider a few special classes of hypersurfaces in Minkowski space $(V,F)$, namely, the \emph{Minkowski
hypersphere} $\hat S_{F}^{n-1}(r):F(x)=r$, the \emph{reverse
Minkowski hypersphere} $\hat S_{F_-}^{n-1}(r):F(-x)=r$ and
\emph{$F^*$-Minkowski cylinder} $\hat
S_{\tilde{F}_\pm}^{m-1}(r)\times \mathbb{R}^{n-m}$, where  $\hat
S_{\tilde{F}_\pm}^{m-1}(r)$ is the (reverse) Minkowski hypersphere
in $(\bar V,\tilde{F})$, $\tilde{F}$ is the dual metric of
$F^*|_{\bar V^*}$ in  $\bar V$,  $\bar V^*=\bar V$ is an
$m$-dimensional subspace of $V^*=V$ and $F^*$ is the dual metric of
$F$. In general, $\tilde{F}\neq F|_{\bar V}$( see Section 5 below
for details).

\begin{thm}\label{thm1-2}In an $n$-dimensional Minkowski space $(V,F)$ with
the BH(resp. HT)-volume form $d\mu$, all the hyperplane, the
Minkowski hypersphere $\hat S_{F_\pm}^{n-1}(r)$ and the
$F^*$-Minkowski cylinder $\hat S_{\tilde{F}_\pm}^{m-1}(r)\times
\mathbb{R}^{n-m}$ must be isoparametric hypersurfaces with one or
two distinct constant principal curvatures.
\end{thm}

\begin{thm}\label{thm1-3}In an $n$-dimensional Randers-Minkowski space $(V,F)$ with
the BH(resp. HT)-volume form $d\mu$,  any isoparametric hypersurface
has two distinct principal curvatures at most, which must be either
a Minkowski hyperplane, a (reverse)  Minkowski hypersphere $\hat
S_{F_{\pm}}^{n-1}(r)$, or an (reverse) $F^*$-Minkowski cylinder
$\hat S_{\tilde{F}_\pm}^{m-1}(r)\times \mathbb{R}^{n-m}$.
\end{thm}

Theorem 1.1 $\thicksim$ 1.3 partially generalize the corresponding
results in \cite{27} and \cite{Ca}. In fact, we also obtain some
generalizations for more general case(see Section 4 and Section 5
for details).

The contents of this paper is organized as follows. In Section 2,
some fundamental concepts and formulas which are necessary for the
present paper are given.  The Gauss-Weingarten formulas and the
Gauss-Codazzi equations of a hypersurface  in $(M,F,d\mu)$ are
established in Section 3. The Finslerian isoparametric hypersurfaces
are discussed in Section 4. In Section 5 and Section 6, we obtain a
Cartan type formula similar to (\ref{1.1})  and give some
classification theorems of isoparametric hypersurfaces in Minkowski
spaces. Particularly,  for the case of Finslerian isoparametric
hypersurfaces in Randers-Minkowski spaces, we give the complete
classification. Finally, in Example 4 of Section 6, we construct a
function which is  transnormal but not isoparametric. It shows that
Theorem B in \cite{WQ} does not hold in Finsler geometry.

\vspace{3mm}

\section{  \textbf{Preliminaries}}

\subsection{ Finsler manifolds}

Throughout this paper, we assume that $M$ is an $n$-dimensional
oriented smooth manifold. Let $TM$ be the tangent bundle over $M$
with local coordinates $(x,y)$, where $x=(x^1,\cdots ,x^{n})$ and
$y=(y^1,\cdots ,y^{n})$. A {\it Finsler metric} on $M$ is a function
$F: TM\longrightarrow[0,\infty)$ satisfying the following
properties: (i) $F$ is smooth on $TM\backslash\{0\}$; (ii)
$F(x,\lambda y)=\lambda F(x,y)$ for all $\lambda>0$; (iii) the
induced quadratic form $g$ is positive-definite, where
$$g:=g_{ij}(x,y)dx^{i} \otimes dx^{j}, ~~~~~~~~~   g_{ij}(x,y)=\frac{1}{2}[F^{2}] _{y^{i}y^{j}}.  $$
Here and from now on, we will use the following convention of index
ranges unless other stated:
$$1\leq i, j,\cdots \leq n ;~~~~1\leq a, b, \cdots \leq n-1.$$
The projection $\pi : TM\longrightarrow M $ gives rise to the
pull-back bundle $\pi^{\ast}TM$ and its dual bundle
$\pi^{\ast}T^{\ast}M$ over $TM\backslash\{0\}$. In
$\pi^{\ast}T^{\ast}M$ there is a global section
$\omega=[F]_{y^{i}}dx^{i},$ called the {\it Hilbert form}, whose
dual is $\ell=\ell^{i}\frac{\partial}{\partial
x^{i}},\ell^{i}=\frac{y^{i}}{F}$, called the {\it distinguished
field}.

As is well known, on the pull-back bundle $\pi^{\ast}TM$ there
exists uniquely the \emph{Chern connection} $\nabla$ with $\nabla
\frac{\partial}{\partial x^i}=\omega_{i}^{j}\frac{\partial}{\partial
x^{j}}$ satisfying (\cite {BCS})
$$\omega^{i}_{j}\wedge dx^{j}=0,$$
$$dg_{ij}-g_{ik}\omega^{k}_{j}-g_{kj}\omega^{k}_{i}=2FC_{ijk}\delta y^{k},\qquad \delta y^{k}:=\frac{1}{F}(dy^k+y^j\omega^{k}_{j}),$$
where $C_{ijk}=\frac{1}{2}\frac{\partial g_{ij}}{\partial y^k}$ is
called the\emph{ Cartan tensor}. The curvature 2-forms of the Chern
connection $\nabla$ are
\begin{eqnarray}\label{2.1}
d\omega^{i}_{j}-\omega^{k}_{j}\wedge
\omega^{i}_{k}=\Omega^{i}_{j}:=\frac{1}{2}R^{i}_{j~kl}dx^{k}\wedge
dx^{l}+P^{i}_{j~kl}dx^{k}\wedge\delta y^{l},
\end{eqnarray}
where $R^{i}_{j~kl}=-R^{i}_{j~lk}$. The {\it flag curvature tensor}
is defined by
$$R^i_{~k}:=\ell^jR^{i}_{j~kl}\ell^l,\qquad R_{jk}=g_{ij}R^i_{~k}.$$
For a unit vector $V=V^{i}\frac{\partial}{\partial x^i}$ with
$g_{ij}y^iV^j=0$, the \emph{flag curvature} $K(y;V)$ is defined by
$$K(y;V)=R_{ij}V^{i}V^{j}.$$
The \emph{Ricci curvature} for $(M,F)$ is defined as
\begin{align}D^{w}_{v}X(x):=\left\{v^{j}\frac{\partial X^{i}}{\partial x^{j}}(x)+\Gamma^{i}_{jk}(w)v^{j}X^{k}(x)\right
\}\frac {\partial}{\partial x^{i}},\label{2.2}\end{align} where
$\{e_1,\cdots
,e_{n}=\ell\}$ is an orthonormal basis with respect to $g_y$.\\

Let $X=X^{i}\frac{\partial}{\partial x^{i}}$ be a vector field. Then
the \emph{covariant derivative} of $X$ along
$v=v^i\frac{\partial}{\partial x^i}$ with respect to $w\in
T_{x}M\backslash \{0\}$ is defined by
\begin{align}D^{w}_{v}X(x):=\left\{v^{j}\frac{\partial X^{i}}{\partial x^{j}}(x)+\Gamma^{i}_{jk}(w)v^{j}X^{k}(x)\right
\}\frac {\partial}{\partial x^{i}},\label{Z1}\end{align} where
$\Gamma^{i}_{jk}$ denotes the connection coefficients of the Chern
connection.

Let ${\mathcal L}:TM\longrightarrow T^{\ast}M$ denote the {\it
Legendre transformation}, which satisfies ${\mathcal L}( 0)=0$ and
${\mathcal L}(\lambda y)=\lambda {\mathcal L}(y)$ for all
$\lambda>0$ and $ y\in TM\setminus \{0\}$. Then ${\mathcal
L}:TM\setminus \{0\}\longrightarrow T^{\ast}M\setminus \{0\}$ is a
norm-preserving $C^\infty$ diffeomorphism  and (\cite{SZ})
\begin{align}\mathcal L(y)=F(y)[F]_{y^{i}}(y)dx^i=F(y)\omega,~~~~~~~~y\in TM\setminus 0.\label{Z01}\end{align}
For a smooth function $f:M\longrightarrow R$, the \emph{gradient
vector} of $f$ at $x\in M$ is defined as $\nabla f(x):={\mathcal
L}^{-1}(df(x))\in T_{x}M$, which can be written as
\begin{align}\nabla f(x):=\left\{\begin{array}{l}
g^{ij}(x,\nabla f)\frac{\partial f}{\partial x^{j}}\frac{\partial }{\partial x^{i}},~~~~~~ df(x) \neq 0,\\
0,~~~~~~~~~~~~~~~~~~~~~~~~~~~~~~~~df(x)=0.
\end{array}\right.\label{Z2}\end{align}
Set $M_{f}:=\{x\in M|df(x)\neq 0\}$.  We define $\nabla^{2}f(x)\in
T^{\ast}_{x}M\otimes T_{x}M$ for $x\in M_{f}$ by using the following
covariant derivative
\begin{align}\nabla^{2}f(v):=D^{\nabla f}_{v}\nabla f(x)\in T_{x}M,~~~~~~~~v\in T_{x}M.\label{Z3}\end{align}
 Set
\begin{align}
D^{2}f(X,Y):=g_{\nabla f}(\nabla^{2}f(X),Y)= g_{\nabla f}(D^{\nabla
f}_{X}(\nabla f),Y). \label{2.3}\end{align} Then we have(\cite{YHS})
\begin{align} g_{\nabla f}(D^{\nabla
f}_{X}(\nabla f),Y)=D^{2}f(X,Y)=D^{2}f(Y,X)= g_{\nabla f}(D^{\nabla
f}_{Y}(\nabla f),X),\label{2.4}\end{align} for any $X,Y\in T_xM.$

Let $d\mu=\sigma(x)dx^{1}\wedge\cdots\wedge dx^{n}$ be an arbitrary
volume form on $(M,F)$.  The \emph{divergence} of a smooth vector
field $V=V^{i}\frac{\partial}{\partial x^{i}}$ on $M$ with respect
to $d\mu$ is defined by
\begin{align} \textmd{div}V :=\frac{1}{\sigma}\frac{\partial}{\partial x^{i}}\left(\sigma V^i\right)=
\frac{\partial V^{i}}{\partial x^{i}}+V^{i}\frac{\partial\ln
\sigma}{\partial x^{i}}. \label{Z4}\end{align} Then the
\emph{Finsler-Laplacian} of $f$ can be defined by
\begin{align}\Delta f:=\textmd{div}(\nabla f)=\frac{1}{\sigma}\frac{\partial}{\partial x^{i}}\left(\sigma g^{ij}(\nabla f)f_j\right).\label{2.6}\end{align}

\subsection{Isometric immersion}

Let $(N,\bar F)$ and $(M,F)$ be Finsler manifolds of dimensions $m(<
n)$ and $n$  respectively. An immersion $\phi:(N,\bar F)\to( {M},
{F})$ is called to be \emph{isometric} if $\bar F(u,v)= {F}(\phi
(u),d\phi (v))$ for $(u,v)\in TN$ with local coordinates
$(u^1,\cdots ,u^{m},v^1,\cdots ,v^{m})$. Here and from now on, we
will use the following convention of index ranges unless other
stated:
$$1\leq \alpha, \beta \cdots \leq m.$$
For an isometric immersion $\phi $, we have
\begin{eqnarray*}
\bar g_{\alpha \beta }(u,v)={g}_{ i  j }(  {x},  {y})\phi ^{ i
}_\alpha \phi ^{ j }_\beta ,
\end{eqnarray*}
where
\begin{eqnarray*} \bar g_{\alpha \beta }=\frac 12[{\bar F}^2]_{v^{\alpha}v^{\beta}},\quad
{x}^{ i }=\phi ^{ i }(u),\quad  {y}^{ i }=\phi ^{ i }_\alpha
v^\alpha ,\quad \phi ^{ i }_\alpha =\frac {\partial \phi ^{ i
}}{\partial u^\alpha }.
\end{eqnarray*}

Set
\begin{eqnarray}\label{2.8}
\aligned &h^{ i }=\phi ^{ i }_{\alpha \beta }v^\alpha v^\beta
-\phi^{ i }_{\alpha}\bar G^{\alpha}+  G^{ i },\\
&h(v):=\frac {h^{ i }(v)}{\bar F^2(v)}\frac {\partial}{\partial
{x}^{ i }}= {\nabla}_{\ell}(d\phi (\ell)),
\endaligned
\end{eqnarray}
where $\phi ^{ i }_{\alpha \beta }=\frac {\partial^2\phi ^{ i
}}{\partial u^\alpha \partial u^\beta},$ $\bar G^{\alpha}$ and $G^{
i }$ are the geodesic coefficients of $(N,\bar F)$ and $(M,F)$
respectively and $h$ is the {\it normal curvature vector field}
(\cite{SZ1}). $(N,\bar F)$ is said to be {\it totally geodesic} if
$h=0$. Let
$$\mathcal{V}(N)=\{\xi\in \phi ^{-1}T^{*} {M}|~\xi (d\phi (X))=0,~~\forall X\in TN\},$$
which is called the {\it normal bundle} of $N$ in $(M,F)$.

\section{ \textbf{Theory of hypersurfaces}}
\subsection{Variation of the induced volume}
Now let $\phi :N\to M$ be an embedded hypersurface (i.e., $m=n-1$)
of Finsler manifold $(M,F)$. For simplicity, we will use $(x,y)\in
TN$, where $x=\phi (u),~y=d\phi v$ for $(u,v)\in TN$. For any $x\in
N$, there exist exactly two unit\emph{ normal vectors}
$\textbf{n}_{\pm}$ such that
$$T_{x}(N)=\{X\in T_{x}(M)|~g_{\textbf{n}_{\pm}}(\textbf{n}_{\pm},X)=0,\quad g_{\textbf{n}_{\pm}}(\textbf{n}_{\pm},
\textbf{n}_{\pm})=1. \}$$ In fact, $\textbf{n}_{\pm}$ are just
${\mathcal L}^{-1}(\pm \nu )$, where $\nu\in\mathcal{V}(N)$ is a
unit 1-form. If $F$ is reversible, then
$\textbf{n}_{-}=-\textbf{n}_{+}$.

Let $\textbf{n}={\mathcal L}^{-1}(\nu )$ be a given normal vector of
$N$ in $(M,F)$, and put $\hat g:=g_{\bf n}$. Let
$d\mu_M=\sigma(x)dx^1\wedge\cdots\wedge dx^{n}$ be an arbitrary
volume form on $(M,F)$. The induced  volume form on $N$ determined
by $d\mu_M$ can be defined by
\begin{align}d\mu_{\textbf{n}}&=\sigma_{\textbf{n}}(u)du^1\wedge
\cdots\wedge du^{n-1}\nonumber\\
&={\sigma}({\phi(u)})\phi^*(i_{\bf n}(dx^1\wedge\cdots\wedge
dx^{n})),~~~ ~~~~~~u\in N,\label{Z21}\end{align} where $i_{\bf n}$
denotes the inner multiplication with respect to $\bf n$.

If we put ${x}=\phi(u), z=(z^{ i }_{\alpha })=(\phi^{ i }_{\alpha
})$, then we have
$$\sigma_{\textbf{n}}(u)=\mathcal{F}({x},z).$$
Consider a smooth variation $\phi_t:N\to M$, $t\in (-\varepsilon
,\varepsilon )$, such that $\phi_0=\phi$. $\phi_t$ induces a family
of volume forms $d\mu_t=\sigma_{\textbf{n}_t}(u)du$ on $\phi_t(N)$
with $d\mu_0=d\mu_{\textbf{n}}$ and $\textbf{n}_0=\textbf{n}$, where
$du=du^1\wedge \cdots\wedge du^{n-1}$. Let $
X=X^i\frac{\partial}{\partial x^i}\in TM$ be the {\it variation
field} of $\{\phi_t\}$ and ${\text {Vol}}(t):=\int_{N}d\mu_t$. As
similar to \cite{SZ1}, we can get that   $$\frac{d
{\text{Vol}}(t)}{dt}\Bigr|_{t=0}=-\int_{N}\mathcal{H}_{d\mu_{\textbf{n}}}(
{X})d\mu_{\textbf{n}},$$ where
\begin{align}
\mathcal{H}_{d\mu_{\textbf{n}}}(
{X})=\frac{1}{\mathcal{F}}\left\{\frac{\partial^{2}\mathcal{F}}{\partial
z^{ i }_{\alpha }
\partial z^{ j }_{\beta }}\phi ^{ j }_{\alpha \beta }-\frac{\partial
\mathcal{F}}{\partial {x}^{ i
}}+\frac{\partial^{2}\mathcal{F}}{\partial  {x}^{j}\partial z^{ i
}_{\alpha }}\phi ^{j}_{\alpha}\right\} {X}^{ i }.
\label{2.10}\end{align} $\mathcal{H}_{d\mu_{\textbf{n}}}$ is called
the \emph{$d\mu_{\textbf{n}}$-mean curvature form} of $\phi$ with
respect to $\bf n$. Define
\begin{align}
{H}_{\textbf{n}}:=\mathcal{H}_{d\mu_{\textbf{n}}}(\textbf{n})
\label{2.11}\end{align} We call ${H}_{\textbf{n}}$ the
\emph{$d\mu_{\textbf{n}}$-mean curvature} of $N$ in $(M,F)$.

\begin{rem}\label{rem3-1} The induced volume form (\ref{Z21}) is determined by the
volume form $d\mu_M$ of $(M,F)$ and the normal vector $\textbf{n}$.
In \cite{SZ1} and \cite{[HS1]}, the induced volume forms (e.g.,
BH-volume form and HT-volume form) are determined by the induced
metric $\bar F$ from $F$, respectively. Hence, in general, the mean
curvature form (\ref{2.10}) is  different from that in \cite{SZ1}
and \cite{[HS1]}. In Riemannian case, they are all same.
\end{rem}

\subsection{Gauss-Weingarten formulas}
 For any tangent vector $y\in T_xN$ at $x\in
N$, there is a unique geodesic $\bar\gamma$ in $(N,\bar F)$ such
that the curve $\gamma (t):=\phi\circ\bar\gamma (t)$ satisfies
$\gamma (0)=x$ and $\dot\gamma (0)=y$. Define
\begin{equation}\label{2.12}\Lambda_{\textbf{n}}(y ):=\hat g(\textbf{n},D^{\dot{\gamma}}_{\dot{\gamma}}
\dot{\gamma}(0)),\qquad y\in T_xN,\end{equation} where $D$ is the
covariant derivative defined by (\ref{Z1}).
$\Lambda_{\textbf{n}}(y)$ is called the \emph{normal curvature} of
$N$ with respect to $\textbf{n}$ (\cite{SZ}). From (14.16) in
\cite{SZ}, it is easy to see that
\begin{equation}\label{2.13}\Lambda_{\textbf{n}}(y )=F^2(y)\hat g(\textbf{n},h(y))=F^2(y)\nu(h(y)),~~~~y\in T_xN,\end{equation}
where $h(y)$ is the normal curvature vector field defined by
(\ref{2.8}). It is obvious that $(N,\bar F)$ is totally geodesic if
and only if $\Lambda_{\textbf{n}}(y )=0, \forall y\in T_xN$. In
general, $\Lambda_{\textbf{n}}(y)$ and $h(y)$ are nonlinear.

  The \emph{Weingarten formula} with respect to $\hat g$ is given by
\begin{equation*}D^{\textbf{n}}_{X}\textbf{n}=\hat\nabla^{\perp}_{X}\textbf{n}-\hat{A}_{\textbf{n}}(X),
\end{equation*}
where $X\in \Gamma(T N)$, $\hat\nabla^{\perp}$ is the \emph{induced
normal connection}, and $\hat{A}_{\textbf{n}}:T_xN\rightarrow T_xN$
is called the\emph{ Weingarten transformation} (or \emph{Shape
operator}). It is obvious that
$$\hat\nabla^{\perp}_{X}\textbf{n}=\hat
g(D^{\textbf{n}}_{X}\textbf{n},\textbf{n})\textbf{n}=\frac{1}{2}\left[X\hat
g(\textbf{n}, \textbf{n})\right]\textbf{n}=0.$$ From \cite{SZ}
(P.222) we know that $\hat{A}_{\textbf{n}}$ is linear and
self-adjoint with respect to $\hat g$. Moreover, we call the
eigenvalues of $\hat{A}_{\textbf{n}}$, $k_1,k_2,\cdots,k_{n-1}$ ,
the \emph{principal curvatures} of $N$ with respect to $\textbf{n}$.
If $\hat{A}_{\textbf{n}}(X)=kX, \forall X\in \Gamma(TN)$, or
equivalently, $k_1=k_2=\cdots=k_{n-1}$, we call $N$ a
\emph{$\hat{g}$-totally umbilic} hypersurface of Finsler manifold
$(M,F)$.

 Define
\begin{equation}\label{2.15}\hat{h}(X,Y):=\hat {g}(\hat{A}_{\textbf{n}}(X),Y)=\hat {g}(\textbf{n},D^{\textbf{n}}_{X}Y),
\quad X,Y\in \Gamma(TN).\end{equation} It is clear that $\hat{h}$ is
bilinear and $\hat{h}(X,Y)=\hat{h}(Y,X)$. Moreover, we define
\begin{equation}\label{2.151}\hat{\nabla}_{X}Y:=(D^{\textbf{n}}_{X}Y)^\top
=D^{\textbf{n}}_{X}Y-\hat{h}(X,Y)\textbf{n}, \quad X,Y\in
\Gamma(TN).\end{equation}
 Then it is easy to prove that $\hat{\nabla}$ is a torsion-free linear connection
on $N$ and satisfies
\begin{eqnarray}\label{2.152}
(\hat{\nabla}_X\hat
g)(Y,Z)=2C_{\textbf{n}}(D^{\textbf{n}}_{X}\textbf{n},Y,Z)=-2C_{\textbf{n}}(\hat
A_{\textbf{n}}(X),Y,Z),
\end{eqnarray}
which shows that $\hat{\nabla}$ is not the Levi-Civita connection on
Riemannian manifold $(N,\hat{g})$. But we also have
 the \emph{Gauss-Weingarten formulas}
\begin{align}\label{2.14}D^{\textbf{n}}_{X}Y&=\hat{\nabla}_{X}Y+\hat{h}(X,Y)\textbf{n},\\
D^{\textbf{n}}_{X}\textbf{n}&=-\hat{A}_{\textbf{n}}(X),~~~~~~~~~~~~~~~~~~~~~~\quad
X,Y\in \Gamma(TN) .\label{2.141}\end{align}

Let $\{e_{i}\}_{i=1}^{n}$ be an orthonormal frame field with respect
to $\hat {g}$ such that $e_{n}=\textbf{n}$. Set(\cite{SZ}\cite{YHS})
\begin{equation}\label{2.16}\hat{H}_{\textbf{n}}:=\sum_{a=1}^{n-1}\hat h(e_{a},e_{a}),\end{equation}
which is independent of the choice of the local frame field
$\{e_{a}\}_{a=1}^{n-1}$.   $\hat{h}$ and $\hat{H}_{\textbf{n}}$ are
called the \emph{$\hat{g}$-second fundamental form} and  the\emph{
$\hat{g}$-mean curvature} of $N$ in $(M,F)$, respectively. It
follows from (\ref{2.15}) and (\ref{2.16}) that
\begin{equation}\label{2.17}\hat{H}_{\textbf{n}}=\sum_{a=1}^{n-1}k_{a},\end{equation}
where $k_a$ are principal curvatures of $N$.

The following lemma gives the relationship between the
$\hat{g}$-second fundamental form and the normal curvatures.

\begin{lem}\label{lem3-2}{\rm(\cite{SZ})} Let $(M,F)$ be a Berwald manifold. Then we have
\begin{equation}\label{2.18}\Lambda_{\textbf{n}}(X)=\hat{h}(X,X)=\hat {g}(\hat
A_{\textbf{n}}(X),X),\end{equation} for any $X\in T_{x}( N)$.
\end{lem}

From the Lemma, we have immediately
$$\hat{H}_{\textbf{n}}=\sum_{a=1}^{n-1}\Lambda_{\textbf{n}}(e_{a}),\qquad
\Lambda_{\textbf{n}}(e_{a})=\hat g(\textbf{n},h(e_{a}))=k_a,$$
 where $\{e_{a}\}_{a=1}^{n-1}$ is an  orthonormal frame field comprised of eigenvectors
of $\hat{A}_{\textbf{n}}$.
\begin{lem}\label{lem3-3} Let $(M,F)$ be a Berwald manifold. Then $(N,\bar F)$ is a totally geodesic hypersurfaces if and only if $k_1=k_2=\cdots=k_{n-1}=0$.
\end{lem}

\subsection{Gauss-Codazzi equations}
 Let the curvature tenser of the
connection $\hat{\nabla}$  be
\begin{eqnarray}\label{12.2}
\hat{R}(X,Y)=\hat{\nabla}_X\hat{\nabla}_Y-\hat{\nabla}_Y\hat{\nabla}_X-\hat{\nabla}_{[X,Y]},~~~~\forall
X,Y\in \Gamma(TN).
\end{eqnarray}
Then we have
\begin{align}\label{12.3}\hat{R}(X,Y)Z&=D^{\textbf{n}}_{X}D^{\textbf{n}}_{Y}Z-D^{\textbf{n}}_YD^{\textbf{n}}_{X}Z-D^{\textbf{n}}_{[X,Y]}Z\nonumber\\
&+\hat{h}_{\textbf{n}}(Y,Z)\hat
A_{\textbf{n}}(X)-\hat{h}_{\textbf{n}}(X,Z)\hat
A_{\textbf{n}}(Y)\nonumber
\\
&-(\hat{\nabla}_X\hat{h}_{\textbf{n}})(Y,Z)\textbf{n}+(\hat{\nabla}_Y\hat{h}_{\textbf{n}})(X,Z)\textbf{n}.
\end{align}
On the other hand, we have
\begin{align*}D^{\textbf{n}}_{\frac{\partial}{\partial x^i}}D^{\textbf{n}}_{\frac{\partial}{\partial x^j}}\frac{\partial}{\partial x^k}
&=D^{\textbf{n}}_{\frac{\partial}{\partial
x^i}}\left(\Gamma^l_{jk}(\textbf{n})\frac{\partial}{\partial
x^l}\right)=
\left(\frac{\partial (\Gamma^s_{jk}(\textbf{n}))}{\partial x^i}+\Gamma^l_{jk}(\textbf{n})\Gamma^s_{il}(\textbf{n})\right)\frac{\partial}{\partial x^s}\\
&=\left\{\frac{\delta \Gamma^s_{jk}}{\delta
x^i}(\textbf{n})+\frac{\partial \Gamma^s_{jk}}{\partial
y^l}(\textbf{n})\left(\frac{\partial \textbf{n}^l}{\partial
x^i}+N^l_i(\textbf{n})\right)+\Gamma^l_{jk}(\textbf{n})\Gamma^s_{il}(\textbf{n})\right\}
\frac{\partial}{\partial x^s}.
\end{align*}
 Then for any $\tilde X,\tilde Y,\tilde Z \in TM$,

 \begin{align}
D^{\textbf{n}}_{\tilde X}D^{\textbf{n}}_{\tilde Y}\tilde Z
&-D^{\textbf{n}}_{\tilde Y}D^{\textbf{n}}_{\tilde X}\tilde
Z-D^{\textbf{n}}_{[\tilde X,\tilde Y]}\tilde Z\nonumber
\\
&=R_{\textbf{n}}(\tilde X,\tilde Y)\tilde Z- P_{\textbf{n}}(\tilde
Y,D^{\textbf{n}}_{\tilde X}\textbf{n})\tilde Z
+P_{\textbf{n}}(\tilde X,D^{\textbf{n}}_{\tilde Y}\textbf{n})\tilde
Z,\label{12.4}\end{align} where $R_{\textbf{n}}$ and $
P_{\textbf{n}}$ are the Chern-curvature tensors of $(M,F)$ defined
by (\ref{2.1}) in $y=\textbf{n}$ and $P(\frac{\partial}{\partial
x^k},\frac{\partial}{\partial x^l})\frac{\partial}{\partial
x^j}=P^i_{jkl}\frac{\partial}{\partial x^i}.$
 It follows from (\ref{12.3}) and (\ref{12.4}) that
 \begin{align}
\hat{R}(X,Y)Z&=R_{\textbf{n}}( X , Y ) Z +P_{\textbf{n}}( Y ,\hat
A_{\textbf{n}}(X)) Z
-P_{\textbf{n}}( X ,\hat A_{\textbf{n}}(Y)) Z \nonumber\\
&+\hat{h}_{\textbf{n}}(Y,Z)\hat
A_{\textbf{n}}(X)-\hat{h}_{\textbf{n}}(X,Z)\hat
A_{\textbf{n}}(Y)\nonumber
\\
&-(\hat{\nabla}_X\hat{h}_{\textbf{n}})(Y,Z)\textbf{n}+(\hat{\nabla}_Y\hat{h}_{\textbf{n}})(X,Z)\textbf{n},
\label{12.5}\end{align} for any $X, Y , Z  \in TN$. Since
$\hat{g}(P_{\textbf{n}}( Y ,\hat A_{\textbf{n}}(X)) Z
,\textbf{n})=L(Y ,\hat A_{\textbf{n}}(X),Z)$£¬where $L$ is the
Landsberg curvature of $(M,F)$, we obtain the \emph{Gauss-Codazzi
equations} as follows.

{\bf Theorem 3.1} For the induced connection $\hat{\nabla}$ on
hypersurfaces, we have
\begin{align} \hat{g}(\hat{R}(X,Y)Z,W)&=\hat{g}(R_{\bf{n}}( X , Y )
Z,W) +\hat{g}(P_{\bf{n}}( Y ,\hat A_{\bf{n}}(X)) Z
-P_{\bf{n}}( X ,\hat A_{\bf{n}}(Y)) Z,W) \nonumber\\
&+\hat{h}_{\bf{n}}(Y,Z)\hat{h}_{\bf{n}}(X,W)-\hat{h}_{\bf{n}}(X,Z)\hat{h}_{\bf{n}}(Y,W)\\
\hat{g}(R_{\bf{n}}( X , Y ) Z,\bf{n})&=(\hat{\nabla}_X\hat{h}_{\bf{n}})(Y,Z)-(\hat{\nabla}_Y\hat{h}_{\bf{n}})(X,Z)\nonumber\\
&+L(X,\hat A_{\bf{n}}(Y),Z)-L(Y ,\hat
A_{\bf{n}}(X),Z).\label{2.51}\end{align}

Particularly,  in Minkowski space, the Gauss-Codazzi equations can
be reduced as
  \begin{align}
&\hat{g}(\hat{R}(X,Y)Z,W)=\hat{h}_{\textbf{n}}(Y,Z)\hat{h}_{\textbf{n}}(X,W)-\hat{h}_{\textbf{n}}(X,Z)\hat{h}_{\textbf{n}}(Y,W),\label{12.6}
\\
&(\hat{\nabla}_X\hat{h}_{\textbf{n}})(Y,Z)-(\hat{\nabla}_Y\hat{h}_{\textbf{n}})(X,Z)=0.
\label{12.7}\end{align}

\subsection{Regular level surfaces of functions}

Let $f:U\subset M\to R$ be a $C^{2}$ function  such that
\begin{equation*}\left\{\begin{aligned}
      &N\cap U=\{x\in U|  f (x)=C\};\\
      & d f  (x)\neq0,\quad  x\in N\cap U,
\end{aligned}\right.\end{equation*}
where $U$ is a neighborhood of some $x_{0}\in N$. Then we have $0=Y
f  ={g}_{\nabla f}(\nabla f  ,Y)$ for any $Y\in \Gamma(T(N\cap U
))$. Therefore, $\frac{\nabla f }{F(\nabla f )}\mid_{N\cap U }$ is a
normal vector of $N\cap U $. Here we have used the condition that $d
f  \neq0$ on $N\cap U $.

By choosing $\textbf{n}=\frac{\nabla  f }{F(\nabla  f )}$ , one can
obtain from (\ref{2.3}), (\ref{2.4}) and (\ref{2.15}) that
\begin{align}
\hat{h}(X,Y)=&{g}_{\nabla f}\left(\hat{A}_{\frac{\nabla f }{F(\nabla
f )}}(X),Y\right)=-{g}_{\nabla f}\left(D^{\nabla  f
}_{X}\frac{\nabla
f}{F(\nabla  f )} ,Y\right)\nonumber\\
=&-\frac{1}{F(\nabla  f )}{g}_{\nabla f}\left(D^{\nabla  f
}_{X}\nabla  f ,Y\right)=-\frac{1}{F(\nabla f )}D^{2} f (X,Y),
\label{4.19}\end{align} where $X,Y\in \Gamma(T N)$.
\begin{lem}\label{lem3-6}  Let
$\{e_{i}\}_{i=1}^{n}$ be a local $g_{\nabla f}$-orthonormal basis
such that $e_{n}={\bf{n}}=\frac{\nabla  f }{F(\nabla  f )}$. Then we
have
\begin{align}F(\nabla  f )\hat{H}_{\bf{n}}=-\sum_{a=1}^{n-1}D^{2} f (e_{a},e_{a})=-\sum_{a=1}^{n-1} f _{aa},
\label{2.20}\end{align} where $f_{aa}=D^{2}f(e_{a},e_{a})$.
\end{lem}
\begin{lem}\label{lem3-7}{\rm(\cite{SZ})} Let $(M,F,d\mu)$ be an $n$-dimensional Finsler
manifold and $f:M\to R$  a smooth function. Then on $M_{f}$ we have
\begin{align}\Delta f=\textmd{tr}_{g_{\nabla f}}(D^{2}f)-S(\nabla f)=\sum_{i}f_{ii}-S(\nabla f).
\label{2.22}\end{align}
\end{lem}

\vspace{6mm}

\section{ \textbf{Isoparametric hypersurfaces in a Finsler manifold}}

Similar to the Riemannian case(\cite{CR},\cite{WQ}), we give tha
following

\begin{defn}\label{defn 4-1}
{\rm Let $f$ be a non-constant continuous function defined on a
Finsler manifold $(M,F)$ and  be smooth in $M\setminus S$, where $S$
is a close set with measure zero in $(M,F,d\mu)$.  Set
$M_{f}:=\{x\in M\setminus S|df(x)\neq 0\}$ and $J=f(M_f)$.  $f$  is
said to be an \emph{isoparametric function} on $(M,F,d\mu)$ if there
is a smooth function $\tilde a (t)$ and a  continuous function
$\tilde b (t)$ defined on $J$ such that
\begin{equation}\label{3.1}\left\{\begin{aligned}
      &F(\nabla f)=\tilde a (f),\\
      &\Delta f=\tilde b (f)
\end{aligned}\right.\end{equation}
hold on $M_f$.  All the regular level surfaces $N_t = f^{-1}(t)$
form an \emph{isoparametric family}, each of which is called an
\emph{isoparametric hypersurface} in $(M,F,d\mu)$.  A function $f$
satisfying only the first equation of (\ref{3.1}) is said to be {\it
transnormal}. A geodesic segment $\gamma(s)$ is called an
\emph{$f$-segment} if $f(\gamma(s))$ is an increasing function of
$s$ and $\dot{\gamma}(s)=\frac{\nabla f}{F(\nabla f)}$ in $M_f$. }
\end{defn}
Note that an $f$-segment is necessarily parametrized by its arc
length.
\begin{lem}\label{lem4-2}Let $f$ be an
isoparametric function on $(M,F,d\mu)$ and $\varphi(t)$ be a
non-constant smooth function defined on $J$ satisfying
$\varphi'(t)\geq0$.  Then $\tilde{f}=\varphi \circ f$ is also
isoparametric on $(M,F,d\mu)$.
\end{lem}
\proof If $\varphi'(t)>0$, $g_{\nabla (\varphi \circ f)}=g_{\nabla
f}$. Then in $M_{\tilde{f}}$, we have
$$ F(\nabla (\varphi \circ f))=\varphi'F(\nabla f),~~~~\Delta(\varphi \circ f)=\varphi''F(\nabla f)^2+\varphi'\Delta f.
$$
It is obvious that $\tilde{f}=\varphi \circ f$  also satisfies
(\ref{3.1}).
\endproof

\subsection{ Transnormal functions and parallel level hypersurfaces}

 First of all, we generalize some results of  transnormal functions in Riemannian geometry.
\begin{lem}\label{lem4-3}Let $f$ be a transnormal function on a connected and forward complete Finsler manifold $(M,F)$.
Then on $M_f$, we have the following results\\
(1) There is a function $\rho$ defined on $(M,F)$ such that
$\nabla\rho=\frac{\nabla f}{F(\nabla f)}$.\\
(2) The integral curves of the gradient vector field $\nabla f$ in $M_f$ are all geodesics($f$-segments).\\
(3) The regular level hypersurfaces $N_t=f^{-1}(t)$ are parallel along the direction of $\nabla f$.\\
(4) $d(x,N_{t_2})=d(N_{t_1},y)=\int^{t_2}_{t_1}\frac{1}{\tilde a
(t)}dt$ for any $[t_1,t_2] \in J$, $x\in N_{t_1},~y\in N_{t_2}$ and
the $f$-segments are the shortest curves among all curves connecting
$N_{t_1}$ and $N_{t_2}$.
\end{lem}

\proof (1) Let $F(\nabla f)=\tilde a (f)$ and define functions
\begin{align}s(t)&=\int^{t}_{t_0}\frac{1}{\tilde a (t)}dt, ~~t_0, t \in J ,\nonumber\\
\rho(x)&=s(f(x)), ~~x\in M.\label{3.31}\end{align} Then
$\nabla\rho=\frac{\nabla f}{F(\nabla f)}$, which shows that
$g_{\nabla f}=g_{\nabla \rho}$ and $g_{\nabla \rho}(\nabla
\rho,\nabla \rho)=1$.

(2) Let $\{e_{i}\}_{i=1}^{n}$ be an orthogonal frame field with
respect to $g_{\nabla \rho}$ such that $e_{n}=\nabla \rho$. From
(\ref{2.4}), we have
$$g_{\nabla \rho}(D^{\nabla\rho}_{\nabla\rho}\nabla \rho,e_{a})=g_{\nabla \rho}(D^{\nabla\rho}_{e_{a}}\nabla \rho,
\nabla \rho)=\frac{1}{2}e_{a}\left[g_{\nabla \rho}(\nabla
\rho,\nabla \rho)\right]=0,$$
$$g_{\nabla \rho}(D^{\nabla\rho}_{\nabla \rho}\nabla \rho,\nabla \rho)=\frac{1}{2}e_{n}\left[g_{\nabla \rho}(\nabla
\rho,\nabla \rho)\right]=0.$$ Thus $D^{\nabla\rho}_{\nabla
\rho}\nabla \rho=0$, that is, all the integral curves of  $\nabla
\rho$ are normal geodesics.

(3) Since $N_t=f^{-1}(t)=\rho^{-1}(s(t))$, the regular level
hypersurfaces $N_t=f^{-1}(t)$ are parallel along the direction of
$\nabla f$.

(4) For any $t_1,t_2 \in J ,~~ t_1<t_2, $ and any $x_1\in
N_{t_1},x_2\in N_{t_2}$, let $\sigma: [0,l]\rightarrow M_f$ be a
piecewise $C^1$ curve with $\sigma(0)=x_1$ and $\sigma(l)=x_2$. Then
\begin{align}L(\sigma)=&\int^{l}_{0}F(\dot{\sigma}(s))ds\geq\int^{l}_{0}g_{\nabla \rho}(\nabla \rho,\dot
{\sigma}(s))ds\nonumber\\
=&\int^{l}_{0}\frac{1}{\tilde a (f(\sigma))}df(\dot{\sigma}(s))ds=\int^{t_2}_{t_1}\frac{1}{\tilde a (t)}dt\nonumber\\
=&s(t_2)-s(t_1).\nonumber\end{align} The equality holds if and only
if $\sigma(s)$ is an $f$-segment.
\endproof

\begin{rem}\label{rem4-4} {\rm In a Finsler manifold, there is a bit difference in the definition of parallel hypersurfaces.
Because in general,  $d(x_1,x_2)\neq d(x_2,x_1)$ unless the Finsler
metric is reversible. That is, $N_1$ is parallel to $N_2$ does not
mean that $N_2$ is parallel to $N_1$. }\end{rem}

From (\ref{2.4}) we see that  $\nabla^2 f$ is self-adjoint with
respect to $g_{\nabla f}$ and
$$ \nabla^2 f(\nabla \rho)=D^{\nabla\rho}_{\nabla \rho}\left(\tilde a (f)\nabla \rho\right)=\tilde a (f)\tilde a '(f)
\nabla \rho,
$$
 which implies taht $\nabla \rho$ is an eigenvector. So we can choose  eigenvectors of $\nabla^2 f$ to form an orthogonal frame
$\{e_{i}\}_{i=1}^{n}$ with respect to $g_{\nabla \rho}$ such that
$e_{n}=\nabla \rho$. By using (\ref{2.141}), we know that
$$ \lambda_ae_a=\nabla^2 f(e_a)=D^{\nabla\rho}_{e_a}(\tilde a (f)\nabla
\rho)=\tilde a (f)D^{\nabla\rho}_{e_a}\nabla \rho=-\tilde a
(f)\hat{A}_{\textbf{n}}(e_a),
$$
which shows that each of $\{e_{a}\}_{a=1}^{n-1}$ is also an
eigenvector of $\hat{A}_{\textbf{n}}$.   Thus we have

\begin{lem}\label{lem4-5}  Let $f$ be a transnormal function on Finsler manifold $(M,F)$.
Then $\nabla f$ is an eigenvector of $\nabla^2 f$. Moreover, if
$\lambda_1,\lambda_2,\ldots,\lambda_{n}$ is the eigenvalues of
$\nabla^2 f$, where $\nabla^2 f(\nabla f)=\lambda_{n}\nabla f$, then
\begin{align}\lambda_a|_{N_t}=-\tilde a (f)k_a,~~\lambda_n=\tilde a
(f)\tilde a '(f),\label{3.2}\end{align} where
$k_1,k_2,\ldots,k_{n-1}$ are the principal curvatures of $N_t$.
\end{lem}

We now prove the following

\begin{lem}\label{lem4-6} Let $f$ be a transnormal function on Finsler manifold $(M,F,d\mu)$, then on $N_{t}$, we have
\begin{align}\Delta f=-\tilde a (t){H}_{\bf{n}}+\tilde a '(t)\tilde a (t)\label{Z31},\end{align}
\begin{align}{H}_{\bf{n}}=\hat{H}_{\bf{n}}+\frac{S(\nabla f)}{F(\nabla  f )},\label{3.3}\end{align}
  where $\bf{n}=\frac{\nabla  f }{F(\nabla  f )}$, ${H}_{\bf{n}}$ and $\hat{H}_{\bf{n}}$ are the $d\mu_{\bf{n}}$-mean curvature and $\hat{g}$-mean curvature  defined by (\ref{2.11}) and (\ref{2.16}),  respectively.
\end{lem}
\proof
  By Lemma 4.2(1), $\textbf{n}=\nabla\rho=\rho^i\frac {\partial}{\partial x^i}$, $\nu=d\rho=\rho_idx^i$. Take a new special local coordinate system $\{(u^a,s)\}$ in a neighborhood $U$ such that $\rho(x(u^a,s))=s,$ that is $d\rho=ds$. Then  $$N_t|_U=\{(u^a,s))|s=s(t)=\int^{t}_{t_0}\frac{1}{\tilde a (t)}dt\}.$$  Set $dx^i=z^i_adu^a+z^i_nds$, $(z^i_j)^{-1}=(w^i_j)$ and let $\phi:N_t\to(M,F)$ be embedding. Then located on $N_t$, we have $$z^i_a=\phi ^{i}_a,~~ z^i_n=\rho^i,~~ w^n_i=\rho_i,~~w^a_i\rho^i=0 $$ and $$d\mu_M=\sigma(x)dx^1\wedge\cdots\wedge dx^{n}=\sigma({x(u^a,s)})\text{det}(z^i_j)du\wedge ds.$$
 From (\ref{Z21}), we have
 \begin{align}d\mu_{\textbf{n}}=&\sigma_{\textbf{n}}(u)du=\sigma({\phi(u)})\text{det}(z^i_j)du,\label{Z32}\end{align}
that is,  $\mathcal{F}(x^i,z^i_a)=\sigma(x)\text{det}(z^i_j)$.
Setting $\eta=\text{det}(z^i_j)$ and using (\ref{2.10}) and
(\ref{2.11}), we have
 \begin{align}
{H}_{\textbf{n}}=\frac{1}{\mathcal{F}}\left\{\frac{\partial^{2}\mathcal{F}}{\partial
z^{i}_{a}\partial
z^{j}_{b}}\phi^{j}_{ab}-\frac{\partial\mathcal{F}}{\partial
{x}^{i}}+\frac{\partial^{2}\mathcal{F}}{\partial {x}^{j}\partial
z^{i}_{a}}\phi^{j}_{a}\right\}\rho^i. \label{Z33}\end{align}
\begin{align*}\frac{\partial\mathcal{F}}{\partial z^{i}_{a}}=&\sigma\frac{\partial\eta }{\partial z^{i}_{a}}=\sigma(x)\eta w^a_i,\\
\frac{\partial\mathcal{F}}{\partial
x^{i}}=&\eta\left(\frac{\partial\sigma}{\partial x^{i}} +\sigma(x)
w^k_l\frac{\partial z^l_k}{\partial x^{i}}\right).
\end{align*}
\begin{align*}
\frac{\partial^{2}\mathcal{F}}{\partial z^{i}_{a}\partial z^{j}_{b}}\phi^{j}_{ab}=&\sigma(x)\eta (w^a_iw^b_j-w^a_jw^b_i)\phi^{j}_{ab}=0,\\
\frac{\partial^{2}\mathcal{F}}{\partial {x}^{j}\partial z^{i}_{a}}\phi^{j}_{a}\rho^i=&\eta w^a_i \left(\frac{\partial\sigma}{\partial x^{j}}+\sigma(x) \rho_k\frac{\partial\rho^k}{\partial x^{j}}\right)\phi^{j}_{a}\rho^i-\sigma\eta w^k_i\frac{\partial z^l_k}{\partial x^{j}}w^a_l\phi^{j}_{a}\rho^i\\
=&-\sigma\eta\left(\rho_i\frac{\partial\rho^l}{\partial
x^{l}}-\rho_j\rho_k\frac{\partial\rho^k}{\partial
x^{i}}\right)\rho^i.
\end{align*}
Thus
\begin{align*}
{H}_{\textbf{n}}=-\frac{1}{\sigma}\left(\frac{\partial\sigma}{\partial
x^{i}}\rho^i+\sigma(x)\frac{\partial\rho^i}{\partial
x^{i}}\right)=-\Delta \rho.
\end{align*}
It follows that
 $$
\Delta f=\textmd{div}(\tilde a (t)\nabla  \rho)=-\tilde a
(t){H}_{\textbf{n}}+\tilde a '(t)\tilde a (t).$$
 On the other hand, from Lemma 3.3, Lemma 3.4 and  Lemma 4.3, we have
 \begin{align}\Delta f=&\textmd{tr}_{g_{\nabla f}}(\nabla^{2}f)-S(\nabla f)\nonumber\\
 =&\sum_{i}\lambda_i-S(\nabla f)\nonumber\\
 =&-\tilde a (t)\hat{H}_{\textbf{n}}+\tilde a (t)\tilde a '(t)-S(\nabla f),\nonumber\end{align}
 It is obvious that (\ref{3.3}) holds.
\endproof

\subsection{ Isoparametric functions and isoparametric hypersurfaces}

From (\ref{3.1}) and Lemma 4.4, we obtain immediately the following

{\bf Theorem 4.1} On a Finsler manifold $(M,F,d\mu)$, a transnormal
function  $f$ is isoparametric if and only if each regular level
hypersurface  $N_t$ of  $f$ has constant $d\mu_{\bf{n}}$-mean
curvature ${H}_{\bf{n}}$, where ${\bf{n}}=\frac{\nabla  f }{F(\nabla
f )}$. Particularly, when $M$ has constant {\bf S}-curvature, a
transnormal function  $f$ is isoparametric if and only if  $N_t$ has
constant $\hat{g}$-mean curvature $\hat{H}_{\bf{n}}$.

Further, we have the following
\begin{lem}\label{lem4-8} Let $f$ be an isoparametric function on Finsler manifold $(M,F,d\mu)$  with constant {\bf S}-curvature $(n+1)cF$,  $\rho$ be a function defined by (\ref{3.31}) and  $e_1,\cdots,e_{n-1},e_n=\nabla\rho$ be the eigenvectors of $\nabla^2f$. Then  we have
 \begin{align}
\sum_{a=1}^{n-1}k_a&=\tilde a '(t)-\frac{\tilde b (t)}{\tilde a (t)}-(n+1)c,\label{3.4}\\
\frac{\partial k_a}{\partial \rho}&=K(\nabla\rho;e_a)+k_a^2,\quad
a=1,\cdots,n-1,\label{3.8}\end{align} where $K(\nabla\rho;e_a)$ is
the flag curvature of $(M,F)$.
\end{lem}
\proof (\ref{3.4}) follows immediately from (\ref{3.1}), (\ref{Z31})
and (\ref{3.3}).

 On the other hand, observe that
$$k_a=g_{\nabla\rho}(\hat
A_{\textbf{n}}(e_a),e_a)=-g_{\nabla\rho}(D^{\nabla\rho}_{e_a}\nabla\rho,e_a),\quad
a=1,\cdots,n-1,$$
  and note that $D^{\nabla\rho}_{\nabla\rho}\nabla\rho=0$. From (\ref{12.4}) and (\ref{2.152}), we have
\begin{align}\label{3.a}
-\frac{\partial k_a}{\partial\rho}=&\nabla\rho\left(g_{\nabla\rho}(D^{\nabla\rho}_{e_a}\nabla\rho,e_a)\right)\nonumber\\
=&g_{\nabla\rho}(D^{\nabla\rho}_{\nabla\rho}D^{\nabla\rho}_{e_a}\nabla\rho,e_a)
+g_{\nabla\rho}(D^{\nabla\rho}_{e_a}\nabla\rho,D^{\nabla\rho}_{\nabla\rho}e_a)
+2C_{\nabla\rho}(D^{\nabla\rho}_{e_a}\nabla\rho,e_a,D^{\nabla\rho}_{\nabla\rho}\nabla\rho)\nonumber\\
=&g_{\nabla\rho}(D^{\nabla\rho}_{e_a}D^{\nabla\rho}_{\nabla\rho}\nabla\rho,e_a)
+g_{\nabla\rho}(D^{\nabla\rho}_{[\nabla\rho,e_a]}\nabla\rho,e_a)\nonumber\\
&+g_{\nabla\rho}(R_{\nabla\rho}(\nabla\rho,e_a)\nabla\rho,e_a)
+g_{\nabla\rho}(D^{\nabla\rho}_{e_a}\nabla\rho,D^{\nabla\rho}_{\nabla\rho}e_a)\nonumber\\
=&-K(\nabla\rho;e_a)
+g_{\nabla\rho}(D^{\nabla\rho}_{[\nabla\rho,e_a]}\nabla\rho,e_a)
+g_{\nabla\rho}(D^{\nabla\rho}_{e_a}\nabla\rho,D^{\nabla\rho}_{\nabla\rho}e_a).
\end{align}
 Moreover, by the torsion freeness of the Chern connection,
it follows from (\ref{3.a}) and (\ref{2.4}) that
\begin{align}
-\frac{\partial
k_a}{\partial\rho}=&-K(\nabla\rho;e_a)+g_{\nabla\rho}(D^{\nabla\rho}_{D^{\nabla\rho}_{\nabla\rho}e_a-D^{\nabla\rho}_{e_a}\nabla\rho
}\nabla\rho,e_a)
+g_{\nabla\rho}(D^{\nabla\rho}_{e_a}\nabla\rho,D^{\nabla\rho}_{\nabla\rho}e_a)\nonumber\\
=&-K(\nabla\rho;e_a)+g_{\nabla\rho}(D^{\nabla\rho}_{
e_a}\nabla\rho,D^{\nabla\rho}_{\nabla\rho}e_a-D^{\nabla\rho}_{e_a}\nabla\rho)
+g_{\nabla\rho}(D^{\nabla\rho}_{e_a}\nabla\rho,D^{\nabla\rho}_{\nabla\rho}e_a)\nonumber\\
=&-K(\nabla\rho;e_a)+2
g_{\nabla\rho}(D^{\nabla\rho}_{e_a}\nabla\rho,D^{\nabla\rho}_{\nabla\rho}e_a)
-g_{\nabla\rho}(D^{\nabla\rho}_{ e_a}\nabla\rho,D^{\nabla\rho}_{
e_a}\nabla\rho).\nonumber
\end{align}
Note that
\begin{align}
2
g_{\nabla\rho}(D^{\nabla\rho}_{e_a}\nabla\rho,D^{\nabla\rho}_{\nabla\rho}e_a)
&=2 g_{\nabla\rho}(k_ae_a,D^{\nabla\rho}_{\nabla\rho}e_a)=k_a\nabla\rho(g_{\nabla\rho}(e_a,e_a))=0,\nonumber\\
g_{\nabla\rho}(D^{\nabla\rho}_{ e_a}\nabla\rho,D^{\nabla\rho}_{
e_a}\nabla\rho) &=g_{\nabla\rho}(k_ae_a,k_ae_a)=k_a^2.\nonumber
\end{align}
Thus, we obtain (\ref{3.8}).
\endproof

{\bf Theorem 4.2} Let $(M,F,d\mu)$ be an $n$-dimensional Finsler
manifold with constant flag curvature and constant {\bf
S}-curvature. Then a transnormal function $f$ is isoparametric if
and only if each regular level surface  of  $f$ has constant
principal curvatures.

\proof By Theorem 4.1, it suffices to prove that each regular level
surface  of an isoparametric function  $f$  has constant principal
curvatures. From (\ref{3.8}), we have
\begin{align}\label{3.b}
\frac{\partial}{\partial
\rho}\sum_{a=1}^{n-1}k_a=\text{Ric}(\nabla\rho)+\sum_{a=1}^{n-1}k_a^2.
\end{align}
Set $K(\nabla\rho,e_a)=C$. We know  from (\ref{3.4}) and
(\ref{3.31}) that $\sum_{a=1}^{n-1}k_a$ is only a function of $f=t$
and  (\ref{3.b}) can be rewritten  as
\begin{align} \label{3.c}
\tilde{a}(t)\frac{\partial}{\partial
t}\sum_{a=1}^{n-1}k_a=\frac{\partial}{\partial
\rho}\sum_{a=1}^{n-1}k_a=(n-1)C+\sum_{a=1}^{n-1}k_a^2.
\end{align}
This indicates that $\sum\limits_{a=1}^{n-1}k_a^2$ is also a
function of $t$. It is easily seen from (\ref{3.8}) that
 \begin{align}
\tilde{a}(t)\frac{\partial}{\partial
t}\sum_{a=1}^{n-1}k_a^2=\frac{\partial}{\partial
\rho}\sum_{a=1}^{n-1}k_a^2=2C\sum_{a=1}^{n-1}k_a+2\sum_{a=1}^{n-1}k_a^3.\nonumber
\end{align}
 This implies that $\sum\limits_{a=1}^{n-1}k_a^3$ is also a function of $f$, and so on. By the properties of symmetric polynomials,
we conclude that $k_a$ is constant on $N_{t}$ for any $a$.
\endproof
 The same result can be obtained by (\ref{3.4}) and (\ref{3.b}) if $n = 3$ and $(M, F)$ has constant Ricci curvature.
{\bf Theorem 4.3} Let $(M,F,d\mu)$ be a $3$-dimensional Finsler
manifold with constant {\bf S}-curvature and constant Ricci
curvature. Then a transnormal function $f$ is isoparametric if and
only if the both principal curvatures of $N_t$ with respect to
$\textbf{n}=\frac{\nabla  f }{F(\nabla  f )}$ are constant.

\section{ \textbf{Isoparametric hypersurfaces in  Minkowski spaces}}

 In this section, we suppose that $(V,F,d\mu)$ is an $n$-dimensional Minkowski space and $d\mu$ is
the BH-volume form or HT-volume form. Then the {\bf S}-curvature of
$(V,F,d\mu)$ vanishes identically. Let $F^*$ be the dual metric of
$F$, which is also a Minkowski  metric, and
$g^{*ij}(\xi)=\frac{1}{2}[F^{*2}(\xi)] _{\xi_{i}\xi_{j}}$. Then
(\ref{3.1}) can be written as

\begin{equation}\label{3.9}\left\{\begin{aligned}
      &F^*(df)=\tilde a (f)\\
      &g^{*ij}(df)f_{ij}=\tilde b (f),
\end{aligned}\right.\end{equation}
where $f_{ij}=\frac {\partial^2f}{\partial x^i\partial x^j}$.

\subsection{ Cartan formulas in Minkowski spaces}
Let $N$ be a hypersurface in $(V,F,d\mu)$, $\{e_{i}\}_{i=1}^{n}$ be
an  orthonormal frame field such that $e_n=\textbf{n}$ and
$e_1,\cdots,e_{n-1}$ be the eigenvectors of $\hat{A}_{\textbf{n}}$.
Set $\hat C_{abc}=C_{\textbf{n}}(e_{a},e_{b},e_{c}) $ and $\hat
\Gamma_{abc}=\hat{g}(\hat{\nabla}_{e_{c}}e_{a},e_{b}).$  It follows
from (\ref{2.152}), (\ref{12.6}) and (\ref{12.7}) that
\begin{align}
&\hat \Gamma_{abc}=2k_c\hat C_{abc}-\hat \Gamma_{bac},~~~~~~~~\forall a,b,c,\label{12.9}\\
&\hat{K}(e_{a}\wedge e_{b})=\hat{g}(\hat{R}(e_{a},e_{b})e_{b},e_{a})=k_ak_b,~~a\neq b,\label{12.91}\\
&0=e_{a}(k_b)\delta_{bc}-k_c\hat \Gamma_{bca}-k_b\hat
\Gamma_{cba}-e_{b}(k_a)\delta_{ac}+k_c\hat \Gamma_{acb}+k_a\hat
\Gamma_{cab}, ~~a\neq b,\label{12.8}\end{align} where
$\hat{K}(e_{a}\wedge e_{b})$ is the sectional curvature of the
hypersurfaces $N$ with respect to the connection $\hat{\nabla}$ and
metric $\hat g$. On the other hand,  we know that
\begin{align}
\hat{K}(e_{a}\wedge e_{b})&=\hat{g}(\hat{\nabla}_{e_{a}}\hat{\nabla}_{e_{b}}e_{b},e_{a})-\hat{g}(\hat{\nabla}_{e_{b}}\hat{\nabla}_{e_{a}}e_{b},e_{a})-\hat{g}(\hat{\nabla}_{[e_{a},e_{b}]}e_{b},e_{a})\nonumber\\
&=e_{a}(\hat \Gamma_{bab})-e_{b}(\hat \Gamma_{baa})+\sum_c\left(\hat
\Gamma_{bcb}\hat \Gamma_{caa}-\hat \Gamma_{bca}\hat
\Gamma_{cab}-\hat \Gamma_{bca}\hat \Gamma_{bac}+\hat
\Gamma_{acb}\hat \Gamma_{bac}\right). \label{12.92}\end{align}
 Let $N$ be an isoparametric hypersurface. Then it follows from  (\ref{12.9}) and (\ref{12.8}) that
\begin{align}
(k_c-k_b)\hat \Gamma_{bca}=(k_c-k_a)\hat \Gamma_{acb},~~~~\forall c.
\label{12.11}\end{align}
 When $k_{\tilde{a}}=k_{a}\neq k_{\tilde{b}}=k_{b}$, one can obtain by (\ref{12.9}) and (\ref{12.11}) that
\begin{align}
\hat \Gamma_{ba\tilde{a}}&=\hat \Gamma_{b\tilde{a}a}=\hat \Gamma_{ab\tilde{b}}=\hat \Gamma_{a\tilde{b}b}=0,\nonumber\\
\hat \Gamma_{ba\tilde{b}}&=2k_b\hat C_{ab\tilde{b}},~~~~\hat \Gamma_{ab\tilde{a}}=2k_a\hat C_{a\tilde{a}b},\label{12.10}\\
\hat \Gamma_{aac}&=k_c\hat C_{aac},~~~~\hat \Gamma_{bbc}=k_c\hat
C_{bbc},~~~~\forall c.\nonumber\end{align} Set
$\mathcal{C}_{ijkl}=\frac{\partial C_{ijk}}{\partial y^l}$ and $\hat
{\mathcal{C}}_{abcd}=\mathcal{C}_{\textbf{n}}(e_{a},e_{b},e_{c},e_{d})$.
We obtain
\begin{align}
e_{a}(\hat \Gamma_{bab})=2k_be_{a}(\hat C_{bab})=2k_b\left(-k_a\hat
{\mathcal{C}}_{aabb}+2\sum_c \hat C_{abc}\hat \Gamma_{bca}+\sum_c
\hat C_{bbc}\hat \Gamma_{aca}\right). \label{12.13}\end{align}
Substituting (\ref{12.13}) and  (\ref{12.91}) into (\ref{12.92}) and
using (\ref{12.10}), we have
\begin{align}
k_ak_b=&-2k_ak_b\hat {\mathcal{C}}_{aabb}+\sum_{k_{\tilde{a}}=k_{a}}\left(2k_b \hat C_{bb\tilde{a}}\hat \Gamma_{a\tilde{a}a}+\hat \Gamma_{b\tilde{a}b}\hat \Gamma_{\tilde{a}aa}\right)\nonumber\\
&+\sum_{k_{\tilde{b}}=k_{b}}\left(4k_b\hat C_{ab\tilde{b}}\hat \Gamma_{b\tilde{b}a}+2k_b \hat C_{bb\tilde{b}}\hat \Gamma_{a\tilde{b}a}-\hat \Gamma_{b\tilde{b}a}\hat \Gamma_{\tilde{b}ab}-\hat \Gamma_{b\tilde{b}a}\hat \Gamma_{ba\tilde{b}}\right)\nonumber\\
&+\sum_{k_c\neq k_a,k_b}\left(4k_b \hat C_{abc}\hat \Gamma_{bca}+2k_b \hat C_{bbc}\hat \Gamma_{aca}\right.\nonumber\\
&\left.+\hat \Gamma_{bcb}\hat \Gamma_{caa}-\hat \Gamma_{bca}\hat
\Gamma_{cab}-\hat \Gamma_{bca}\hat \Gamma_{bac}+\hat
\Gamma_{acb}\hat \Gamma_{bac}\right) , ~~k_a\neq k_b .\nonumber
\end{align}It follows from (\ref{12.10}) and (\ref{12.11}) that
\begin{align}
k_ak_b=&k_ak_b\left(-2\hat {\mathcal{C}}_{aabb}+4\sum_c \hat C_{aac}\hat C_{bbc}\right)\nonumber\\
&+\sum_{k_c\neq k_a,k_b}\left(\frac{k_c-k_a}{k_c-k_b}-\frac{k_a-k_c}{k_a-k_b}-\frac{(k_c-k_a)^2}{(k_c-k_b)(k_a-k_b)}\right)\hat \Gamma^2_{acb}\nonumber\\
&+\sum_{k_c\neq k_a,k_b}k_b\hat C_{abc}\hat \Gamma_{acb}\left(\frac{k_c-k_a}{k_c-k_b}+\frac{k_a-k_c}{k_a-k_b}+\frac{(k_c-k_a)^2}{(k_c-k_b)(k_a-k_b)}\right)\nonumber\\
=&k_ak_b\left(-2\hat {\mathcal{C}}_{aabb}+4\sum_c \hat C_{aac}\hat
C_{bbc}\right)+2\sum_{k_c\neq k_a,k_b}\frac{k_c-k_a}{k_c-k_b}\hat
\Gamma^2_{acb}, ~~k_a\neq k_b . \label{12.171}\end{align}

For any $y,X,Y\in V$, we define the \emph{Cartan curvature} as
\begin{align}
Q_{y}(X,Y)=\frac{2F^2(y)}{g_{y}(X,X)g_{y}(Y,Y)}\left(2\sum_i
C_y(X,X,e_i)C_y(Y,Y,e_i)-\mathcal{C}_y(X,X,Y,Y)\right),
\label{12.14}\end{align} where $\{e_{i}\}_{i=1}^{n}$ is an
orthonormal frame of $V$ with respect to $g_y$. By (\ref{12.171}) we
have
\begin{align}
k_ak_b(1-Q_{\textbf{n}}(e_{a},e_{b}))&=2\sum_{k_c\neq k_a,k_b}\frac{k_c-k_a}{k_c-k_b}\hat \Gamma^2_{acb}, ~~~~~~k_a\neq k_b,~~ g>2,\label{12.151}\\
k_ak_b(1-Q_{\textbf{n}}(e_{a},e_{b}))&=0,
~~~~~~~~~~~~~~~~~~~~~~~~~~~~~~~~~~~~~~~k_a\neq k_b,~~ g=2.
\label{12.15}\end{align}

Here and from now on, we suppose that the isoparametric hypersurface
$N_t$ has $g$ distinct constant principal
 curvatures $\kappa_1,\kappa_2,\cdots,\kappa_{g}$ and the  multiplicity of $\kappa_r$ is $m_r$, $r=1,2,\cdots,{g}$.
\begin{lem}\label{lem5-1} Let $(V,F,d\mu)$ be an $n$-dimensional Minkowski space satisfying $Q_{y}(X,Y)=q(y)\neq1$ and $C_y(X,Y,Z)=0$ for any $y,X,Y,Z\in V$, where $y,X,Y,Z $ are all orthogonal to  each other with respect to $g_y$. Then for any isoparametric hypersurface with principal curvatures $\kappa_1<\kappa_2<\ldots<\kappa_g (g\geq2)$, we have the following\textbf{ Cartan type formula}
\begin{align} \sum_{r\neq s}m_r\frac{\kappa_s\kappa_r}{\kappa_s-\kappa_r}=0,~~~~s=1,\ldots,g.
\label{1.11}\end{align}
\end{lem}
\proof It is obvious when $g=2$.  When $g>2$, we know from
(\ref{12.9}) that $ \Gamma_{abc}= -\Gamma_{bac}$, for any $a\neq
b\neq c\neq a$, and (\ref{12.151}) can be written as
\begin{align}
k_ak_b(1-q(\textbf{n}))=2\sum_{k_c\neq
k_a,k_b}\frac{k_c-k_a}{k_c-k_b}\hat \Gamma^2_{acb}, ~~~k_a\neq k_b.
\label{12.16}\end{align} By (\ref{12.9}) and (\ref{12.11}), we have
\begin{align}
(k_a-k_c)\hat \Gamma_{acb}=(k_c-k_b)\hat \Gamma_{cba}=(k_b-k_a)\hat
\Gamma_{bac}. \label{12.17}\end{align}
 Denote $$\lambda_{abc}:=\frac{(k_a-k_b)}{\sqrt{|1-q(\textbf{n})|}}\hat \Gamma_{abc},$$
$$\rho_{rst}:= \sum_{k_a=\kappa_r, k_b=\kappa_s,k_c=\kappa_t}\frac{\lambda^2_{abc}}{(k_b-k_a)(k_a-k_c)(k_b-k_c)},
$$
where $1\leq r,s,t\leq g$ and $r,s,t$ are distinct. Then
$\rho_{rst}$ are skew-symmetric with respect to any two indexes.
Thus (\ref{12.15}) can be written as
\begin{align}
\frac{m_rm_s\kappa_r\kappa_s}{\kappa_s-\kappa_r}=\pm\sum_{t\neq
r,s}\rho_{rst}:=\rho_{rs}, ~~~~r\neq s. \label{12.18}\end{align}
 It follows that
 \begin{align} m_s\sum_{r\neq s}m_r\frac{\kappa_s\kappa_r}{\kappa_s-\kappa_r}=\pm\sum_{r\neq s}\rho_{rs}=\pm\sum_{r,t\neq s}\rho_{rst}=0,~~~~s=1,\ldots,g.
\label{12.19}\end{align}
\endproof
\begin{rem}\label{rem5-2}
There are non-Euclidean Minkowski metrics, for example
Randers-Minkowski metrics,  satisfying the conditions in Lemma
5.1(see Section 6 for details).
\end{rem}
 When $g>2$, using (\ref{12.18}) and (\ref{12.19}), a similar argument as
 in \cite{Ca} yields  a contradiction. Thus we have

{\bf Theorem 5.1}  Let $(V,F,d\mu)$ be an $n$-dimensional Minkowski
space satisfying $Q_{y}(X,Y)=q(y)\neq1$ and $C(X,Y,Z)=0$ for any
$y,X,Y,Z\in V$, where  $y,X,Y,Z $ are all orthogonal to each other
with respect to $g_y$ for $n>3$. Then any isoparametric hypersurface
in $(V,F,d\mu)$ has two distinct principal curvatures at most.

\subsection{Isoparametric hypersurfaces with \textbf{$g=1$} }
In $(V,F)$, the hypersurfaces $\hat S_{F_+}^{n-1}(x_0,r)=\{x\in
V|F(x-x_0)=r\}$ and $\hat S_{F_-}^{n-1}(x_0,r)=\{x\in
V|F(x_0-x)=r\}$ are said to be the \emph{Minkowski hypersphere} and
the \emph{reverse Minkowski hypersphere} of radius $r$ centered at
point $x_0$, respectively.
\begin{rem}\label{rem5-3}If $F$ is reversible, then $\hat S_{F_+}^{n-1}(x_0,r)=\hat S_{F_-}^{n-1}(x_0,r)$. In general, they are different.
\end{rem}
\begin{thm}\label{thm5-4} In an $n$-dimensional Minkowski space $(V,F,d\mu)$,  any
hyperplane,  Minkowski hypersphere and reverse Minkowski hypersphere
with radius $r$ are all isoparametric, which are also
$\hat{g}$-totally umbilic  and have constant principal curvatures
$0$, $\frac{-1}{r}$, $\frac{1}{r}$ and constant sectional curvatures
$0$, $\frac{1}{r^2}$, $ \frac{1}{r^2}$,  respectively.
\end{thm}

\proof For a linear function  $f$ in space $V$,  $F^*(df)$ is a
constant. It is obvious that $f$ satisfies equations (\ref{3.9})
with $\tilde a (t)=$const., $\tilde b (t)=0$ and $\nabla^2 f=0$.
From Lemma 4.3 and  (\ref{12.91}), we see that
$k_1=k_2=\ldots=k_{n-1}=0$ and $\hat K=0$.

Since the Minkowski metric is independent of the point $x\in V$,
that is, $F(x,y)=F(y), y\in T_xV=V$,  it suffices to consider
Minkowski hyperspheres centered at the origin and put
$f(x)=\pm\frac{1}{2}F^2(\pm x), ~\forall x\in V$. Then we have from
(\ref{Z01}) that$$df=FF_{y^i}dx^i|_{y=\pm x}={\mathcal L}(\pm x
),\quad f_{ij}=\pm[\frac{1}{2}F^2]_{y^iy^j}|_{y=\pm x }=\pm
g_{ij}(\pm x ).$$
 Thus $$\nabla  f=\pm x  ,\qquad F^*(df)=F(\pm x  )=\sqrt{\pm 2f(x )},$$
$$g^{*ij}(df)f_{ij}=\pm g^{ij}(\pm x )g_{ij}(\pm x )=\pm n.
$$
The function $f$ also satisfies equations in (\ref{3.9}) with
$\tilde a (f)=\sqrt{\pm2f}$ and $\tilde b (f)=\pm n$. Consider the
isoparametric family $f=t$. By (\ref{3.4}) and (\ref{3.c}), a
straightforward computation yields that
  \begin{align}\label{5.31}\sum_a k_a=\pm\frac{1-n}{\tilde a}
,~~~~~\sum_a k^2_a=\frac{n-1}{\tilde a^2}.
\end{align}
The formulas  above imply that
  \begin{align}\label{5.32}
(\sum_a k_a)^2&=(n-1)\sum_a k^2_a.
\end{align}
From (\ref{5.31}), (\ref{5.32}) and  (\ref{12.91}),  we have
 $$k_1=k_2=\ldots=k_{n-1}=\pm\frac{-1}{\sqrt{\pm2t}}=\mp\frac{1}{{r}},~~~~\hat K=\frac{1}{{r^2}}.
 $$
\endproof
\begin{thm}\label{thm5-5} In an $n$-dimensional Minkowski space $(V,F,d\mu)(n>2)$,  a
$\hat{g}$-totally umbilic hypersurface is isoparametric with $g=1$
and it must be either a Minkowski hyperplane, a Minkowski
hypersphere or a reverse Minkowski hypersphere.
\end{thm}
\proof Let $N$ be a $\hat{g}$-totally umbilic hypersurface  with
local coordinates $(u^1,\cdots ,u^{n-1})$ at $x=\phi(u)\in N$. Note
that $k_1=k_2=\ldots=k_{n-1}=k$ and
 \begin{align}\label{5.25}D^{\textbf{n}}_{\frac{\partial}{\partial u^a}}\textbf{n}=-\hat A_{\textbf{n}}\left(\frac{\partial}{\partial  u^a}\right)=-kd\phi\left(\frac{\partial}{\partial  u^a}\right).\end{align}
  Set $\textbf{n}=n^i\frac{\partial}{\partial x^i}$, $ n^i_a=\frac{\partial n^i}{\partial  u^a} $, $ n^i_{ab}=\frac{\partial^2 n^i}{\partial  u^au^b} $.
Then from  (\ref{5.25}), we obtain
 \begin{align*}n^i_{ab}+\frac{\partial k}{\partial  u^b}\phi^i_{a}+k\phi^i_{ab}=0.\end{align*}
  This implies $\frac{\partial k}{\partial  u^a}=\frac{\partial k}{\partial  u^b}=0, \forall a\neq b$, that is,  $k$ is a constant. Hence, $N$ is isoparametric.

 When $k=0$, (\ref{5.25}) implies that  $\textbf{n}$ is a constant vector. Thus $N$ is a Minkowski hyperplane. When $k\neq0$, from  (\ref{5.25})  we know that there is a point $x_0\in M$ such that
 $$\textbf{n}+kx=kx_0, ~~~\forall x \in N.
 $$
This implies $F(-k(x-x_0))=1$. Therefore, $F(-(x-x_0))=\frac{1}{k}$
when $k>0$ and $F(x-x_0)=\frac{-1}{k}$ when $k<0$.
\endproof
\vspace{4mm}

\subsection{Isoparametric hypersurfaces with $g=2$ }
  Let $(\bar V,\bar F)$ be an $m$-dimensional Minkowski subspace of $(V,F)$, that is, $\bar F=F|_{\bar V}$. Let $\hat S_{\bar{F}_\pm}^{m-1}(\bar x_0,r)$ be the (reverse) Minkowski hypersphere in the $(\bar V,\bar F)$.  We call $\hat S_{\bar{F}_\pm}^{m-1}(\bar x_0,r)\times
\mathbb{R}^{n-m}:\bar{F}(\pm(\bar x-\bar x_0))=r, ~\bar x\in \bar V,
$ the  \emph{(reverse) $F$-Minkowski cylinder} of radius $r$ in
Minkowski space $(V,F)$.

 On the other hand, let  $F^*$ be the dual metric of $F$, $\bar V^*$ be an $m$-dimensional subspace  of $V^*=T^*_x V$ and $\tilde{F}$ be the dual metric of $F^*|_{\bar V^*}$ in $\bar V$. That is,
$$\tilde{F}(\bar{y})=\sup_{\bar{\xi}\in\bar{V}^* \backslash0}\frac{\bar{\xi}(\bar{y})}{{F}^*(\bar{\xi})},
\quad \bar{y}\in\bar{V}.$$ Then $(\bar{V}, \tilde{F})$ is also a
Minkowski space. Let $\hat S_{\tilde{F}_\pm}^{m-1}(\bar x_0,r)$ be
the (reverse) Minkowski hypersphere in  $(\bar V,\tilde{F})$. The
cylinder $\hat S_{\tilde{F}_\pm}^{m-1}(\bar x_0,r)\times
\mathbb{R}^{n-m}:\tilde{F}(\pm(\bar x-\bar x_0))=r, ~\bar x\in \bar
V, $  is said to be the \emph{(reverse) $F^*$-Minkowski cylinder} of
radius $r$ in Minkowski space $(V,F)$.

\begin{rem}\label{rem5-6} In general, we have
 $$\bar F(\bar{y})=F(\bar{y})=\sup_{{\xi}\in{V}^* \backslash0}\frac{{\xi}(\bar{y})}{{F}^*({\xi})}\geq\sup_{\bar{\xi}\in\bar{V}^* \backslash0}\frac{\bar{\xi}(\bar{y})}{{F}^*(\bar{\xi})}
=\tilde{F}(\bar{y}), \quad \bar{y}\in\bar{V}.$$ Thus an (reverse)
$F^*$-Minkowski cylinder is not always an (reverse) $F$-Minkowski
cylinder.
\end{rem}
\begin{thm}\label{thm5-7}In an $n$-dimensional Minkowski space $(V,F,d\mu)$, any (reverse) $F^*$-Minkowski cylinder $\hat S_{\tilde{F}_\pm}^{m-1}(r)\times \mathbb{R}^{n-m}$  must be an isoparametric hypersurfaces with
constant principal curvatures $0$ and $\mp\frac{1}{r}$.
\end{thm}
\proof Let $\bar{V}^* $  be an $m$-dimensional subspace of $V^*$ and
$\tilde{F}^*=F^*|_{\bar{V}^* }$. Then we can take an orthogonal
coordinate system  $\{(x^i)\}$ in $V$ with respect to the Euclidean
metric such that  $\bar{V}^* =
\{\bar{\xi}=(\xi_1,\cdots,\xi_m,0,\cdots,0)|\bar{\xi}\in V^*\}$ and
 $ \bar V=\{x\in V|x=(x^1,\ldots,x^m,0,\ldots,0) \}$.  Here and from now on, we use the following convention of index
ranges:
$$1\leq \alpha, \beta \cdots \leq m<n;~~~m+1\leq \lambda, \mu \cdots \leq n.$$
Set $ \bar x=(x^\alpha,0)\in  \bar V$ and set
$$f(x)=\pm\frac{1}{2}\tilde F^2(\pm\bar x),$$
for any $ x=(x^i)\in V$. Then
\begin{align}
df&=\tilde{F}( \pm\bar x )\tilde{F}_{y^\alpha}(\pm\bar x )dx^\alpha=\tilde{\mathcal{L}}(\pm\bar x)\in\bar{ V}^*,\nonumber\\
f_{ij}(x)&=[\frac{1}{2}\tilde{F}^2]_{y^\alpha y^\beta}(\pm\bar x
)\delta^{\alpha}_i\delta^{\beta}_j
=\pm\tilde{g}_{\alpha\beta}(\pm\bar
x)\delta^{\alpha}_i\delta^{\beta}_j,\nonumber
\end{align}
where $\tilde{\mathcal{L}}:\bar{ V}\to\bar{ V}^*$ is the Legendre
transform  with respect to the metric $\tilde F$. So we have
$$\widetilde{\nabla}  f=\pm\bar x,~~~~F^*(df)=\tilde{F}^*(df)=\tilde{F}(\pm\bar x )=\sqrt{\pm2f(x)},$$
$$g^{*ij}(df)f_{ij}=\pm\tilde{g}^{*\alpha\beta}(df) \tilde{g}_{\alpha\beta}(\pm\bar x)
=\pm\tilde{g}^{\alpha\beta}(\pm\bar x )\tilde
g_{\alpha\beta}(\pm\bar x )=\pm m,
$$
where $\widetilde{\nabla}  f$ is the gradient vector of $f$  in
$(\bar V,\tilde F)$. This means that $f$ satisfies (\ref{3.9}) with
$\tilde a (t)=\sqrt{\pm2t}$, $\tilde b (t)=\pm m$.

 By (\ref{3.8}) and (\ref{3.4}),  a straightforward computation yields that
\begin{align}\label{5.33a}
\sum_a k_a=\pm \frac{1-m}{\tilde a},~~~~\sum_a
k^2_a=\frac{m-1}{\tilde a^2}.
\end{align}
 We get from above that
  \begin{align}\label{5.33b}
(\sum_a k_a)^2&=(m-1)\sum_a k^2_a.
\end{align}
On the other hand, $\nabla^2 f$ has $n-m$ vanishing eigenvalues
since $f_{i\lambda}=f_{\lambda i}=0$. Setting $
k_{m}=k_{m+1}=\ldots=k_{n-1}=0 $, we obtain by (\ref{5.33b}) and
(\ref{5.33a}) that
 $$k_1=k_2=\ldots=k_{m-1}=\pm\frac{-1}{\sqrt{\pm2t}}=\mp\frac{1}{r}.
 $$
\endproof
In general, it is difficult to express the metric $\tilde{F}$
explicitly.  So we consider some special Minkowski spaces in which
some $F$-Minkowski cylinders are also isoparametric.

{\bf Corollary 5.1}  Let $(V,F,d\mu)$ be an $n$-dimensional Minkowski space and $(\bar V,\bar F)$ be a Minkowski subspace of $(V,F,d\mu)$.  If \\
(1) the Legendre transformation satisfies ${\mathcal L}(\bar V)\subset \bar V^*$ or\\
(2) $\bar F$ is a Euclidean metric, \\
then the (reverse) $F$-Minkowski cylinder $\hat S_{\bar
F_\pm}^{m-1}(r)\times \mathbb{R}^{n-m}$ is  isoparametric in
$(V,F,d\mu)$,  where  $\hat S_{\bar F_\pm}^{m-1}(r)$ is a (reverse)
Minkowski hypersphere in  $(\bar V,\bar F)$. In the latter case,
$\hat S_{\bar F_\pm}^{m-1}(r)\times \mathbb{R}^{n-m}$ is actually
the Euclidean cylinder $\mathbb S^{m-1}(r)\times \mathbb{R}^{n-m}$.

\proof

 (1)~~If ${\mathcal L}(\bar V)\subset \bar V^*$ , then

\begin{align*}\bar{F}(\bar{y})&={F}(\bar{y})=\sup_{{\xi}\in{V}^* \backslash0}\frac{{\xi}(\bar{y})}{{F}^*({\xi})}\\
&\geq\sup_{\bar{\xi}\in\bar{V}^*
\backslash0}\frac{\bar{\xi}(\bar{y})}{{F}^*(\bar{\xi})}
=\tilde{F}(\bar{y})\\
&\geq \frac{\bar{\xi}(\bar{y})}
{{F}^*(\bar{\xi})}|_{\bar\xi=\mathcal L(\bar y)}={F}^*(\mathcal
L(\bar y))={F}(\bar{y}), \quad \forall\bar{y}\in\bar{V}.
\end{align*}
That is, $\bar{F}=\tilde{F}$.\\
(2) ~~Let  $\bar F$ be a Euclidean metric. First, for any unit
vector $e$ in Euclidean $m$-space  $(\bar V,\bar F)$, there exists
an orthogonal matrix $\bar A$ such that $\bar
Ae=(\frac{1}{\sqrt{m-1}},\cdots,\frac{1}{\sqrt{m-1}})^T$. In fact,
in the standard sphere $S^m(1)$, by a suitable rotation,  one can
transplant $e$ to
$(\frac{1}{\sqrt{m-1}},\cdots,\frac{1}{\sqrt{m-1}})^T$. Let $\bar
V=\{x\in V|x=(x^1,\ldots,x^m,0,\ldots,0) \}$. Then
$F(y)=F(y^1,\cdots,y^n)$ is invariant under the action of the group
$O(m)\times I(n-m)$ and thus there exists an orthogonal matrix $A\in
O(m)\times I(n-m)$ such that
\begin{align}
&A\left(\frac{y^1}{\sqrt{(y^1)^2+\cdots+(y^m)^2}},\cdots,\frac{y^m}{\sqrt{(y^1)^2+\cdots+(y^m)^2}},y^{m+1}
\cdots,y^n\right)^T\nonumber\\
=&\left(\frac{1}{\sqrt{m}},\cdots,\frac{1}{\sqrt{m}},y^{m+1}\cdots,y^n\right)^T.\nonumber
\end{align}
Thus,
\begin{align}\label{5.2z}
F(y^1,\cdots,y^n)&=F(A(y^1,\cdots,y^m,y^{m+1},\cdots,y^n))\nonumber\\
&=F\left(\frac{\sqrt{(y^1)^2+\cdots+(y^m)^2}}{\sqrt{m}},\cdots,\frac{\sqrt{(y^1)^2+\cdots+(y^m)^2}}{\sqrt{m}},
y^{m+1},\cdots,y^n\right)\nonumber\\
&:=\Phi(\sqrt{(y^1)^2+\cdots+(y^m)^2},y^{m+1},\cdots,y^n).
\end{align}
Let $F^{\ast}(\xi )=F^{\ast}(\xi _1,\cdots,\xi _n)$ be the dual
metric of $F$, where $\xi =\mathcal L(y)$. Define
\begin{align}
\begin{array}{cc}
       A^{\ast}:\quad  &T^{\ast}M \to T^{\ast}M \nonumber\\
         &~~~\xi\longmapsto A^{\ast}\xi
       \end{array}
\end{align}
such that $A^{\ast}\xi (y)=\xi (Ay)$ for any $ y\in TM$. Then
$A^{\ast}=A^T\in O(m)\times I(n-m)$ and we have
$$F^{\ast}(A^{\ast}\xi )=\sup_{y\in TM\backslash0}\frac{A^{\ast}\xi (y)}{F(y)}
=\sup_{y\in TM\backslash0}\frac{\xi (Ay)}{F(y)} =\sup_{Ay\in
TM\backslash0}\frac{\xi (Ay)}{F(Ay)} =F^{\ast}(\xi ).$$ By a similar
argument as above, we obtain
\begin{align}\label{5.1}
F^{\ast}(\xi _1,\cdots,\xi _n):=\Psi(\sqrt{(\xi _1)^2+\cdots+(\xi
_m)^2},\xi _{m+1},\cdots,\xi _n).
\end{align}
Let $\bar{V}^* := \{\bar{\xi}=(\xi_1,\cdots,\xi_m,0,\cdots,0)\}$  be
the subspace of $V^*$. Then $F^*|_{\bar{V}^*}$ is a standard
Euclidean metric and its dual metric $\tilde{F}$ in $\bar V$ is also
a standard Euclidean metric. Thus $\bar F=\tilde F$.
\endproof
It can be seen from (\ref{Z01}) that ${\mathcal L}(\bar V)\subset
\bar V$ if and only if
 \begin{align}[F]_{y^{\lambda}}(\bar y)=0,~~~~ m+1\leq \lambda\leq n,~~~~\forall \bar y\in \bar V.
 \label{4z}\end{align}

{\bf Example 1}  Let $(V,F,d\mu)$ be a Minkowski space with $k$-th
root metric $ F=\sqrt[k]{\sum_i (y^i)^k}, (k>2)$. Set $\bar V=\{x\in
V|x=(x^1,\ldots,x^m,0,\ldots,0) \}(1< m<n)$.  It is easy to check
that
 $F$ satisfies (\ref{4z}), which means that the condition (1) in Corollary 5.1 is satisfied.
Therefore
 $$(x^1)^k+(x^2)^k+\ldots+(x^m)^k=r^k$$  is an isoparametric hypersurface, where $1<m\leq n$.

{\bf Example 2}
   Let $(V,\alpha\phi(\frac{\beta}{\alpha}),d\mu)$ be an $(\alpha,\beta)$-Minkowski space,
   where $\alpha=\sqrt{\sum_i (y^i)^2},\beta=by^1$.
Set $\bar V=\{x\in V|x=(x^1,\ldots,x^m,0,\ldots,0) \}(1< m<n)$. Then
$F$ also satisfies (\ref{4z}), that is, $F$ satisfies the condition
(1) in Corollary 5.1.  On the other hand, if we set $\bar V=\{x\in
V|x=(0,\cdots,0,x^{m+1},\ldots,x^n) \}(1\leq m<n-1)$, then it is
obvious that $\bar F=\sqrt{\sum_\lambda (y^\lambda)^2}$ satisfies
the condition  (2) in Corollary 5.1. Thus
 $$\sqrt{\sum_{\alpha=1}^{m}(x^\alpha)^2}
   \phi\left(\frac{\pm bx^1}{\sqrt{\sum_{\alpha=1}^{m}(x^\alpha)^2} }\right)=r,~~~~1<m\leq n,$$ and
   $$\sum_{\lambda=l+1}^{n}(x^\lambda)^2 =r^2,~~~~1\leq l<n-1,$$
are  isoparametric hypersurfaces in
$(V,\alpha\phi(\frac{\beta}{\alpha}),d\mu)$.

 \begin{thm}\label{thm5-9} Let $(V,F,d\mu)$ be an $n$-dimensional Minkowski space satisfying $Q_{y}(X,Y)\neq1$, $\forall y,X,Y\in V$, where $y,X,Y $ are orthogonal to each other with respect to $g_y$.
Then any isoparametric hypersurface with $g=2$ in $(V,F,d\mu)$ must
be an $F^*$-Minkowski cylinder or a reverse $F^*$-Minkowski
cylinder.
\end{thm}

\proof Let $N$ be an isoparametric hypersurface with $g=2$ and
$\textbf{n}=\nabla\rho$, where $\rho(x)=d(N,x)$. By Lemma 4.1 and
Lemma 4.2, $\rho(x)$ is an isoparametric function in  $(V,F,d\mu)$.
Hence, it satisfies (\ref{3.9}). Let $e_1,\cdots,e_{m-1}$ and
$e_m,\cdots,e_{n-1}$  be the eigenvectors of $\hat A_\textbf{n}$
with constant principal curvatures $k_1=k_2=\ldots=k_{m-1}=k$ and
 $k_{m}=k_{m+1}=\ldots=k_{n-1}=k'$, respectively. By assumption and (\ref{12.151}), $kk'=0$. Set $k'=0, k=\pm\frac{1}{r},r>0$. From (\ref{12.10}), we know that $\hat \Gamma_{\lambda \alpha\beta}=\hat \Gamma_{\lambda \alpha \mu}=0,1\leq\alpha,\beta\leq n-1, m\leq\lambda,\mu\leq n-1$. Thus from (\ref{2.141}), we have
\begin{align*}
D^\textbf{n}_{e_b}e_\lambda&=\sum_{a=1}^{n-1}\hat \Gamma_{\lambda
ab}e_a+\hat h_\textbf{n}(e_b,e_\lambda) =\sum^{n-1}_{\nu=m}\hat
\Gamma_{\lambda \nu b}e_\nu,~~1\leq b\leq n-1,m\leq\lambda\leq n-1.
\end{align*}
 This implies that the tangent subspace spanned by $\{e_{m},\cdots,e_{n-1}\}$ is parallel in $V$. We can
 choose the coordinates $\{x^i\}$ in $V$ such that
$$\text{span}\{e_{m},\cdots,e_{n-1}\}=\text{span}\left\{\frac{\partial}{\partial x^{m+1}},\cdots,\frac{\partial}{\partial x^{n}}\right\}.$$
Then
$$\frac{\partial \rho}{\partial x^\lambda}=d\rho(\frac{\partial}{\partial x^\lambda})=g_{\textbf{n}}(\frac{\partial}{\partial x^\lambda},\textbf{n})=0,\forall \lambda=m+1,\ldots,n,$$
which shows that $\rho(x)=\rho(x^1,\ldots,x^m)$. Set $ \bar V=\{\bar
x\in V|\bar x=(x^1,\ldots,x^m,0,\ldots,0) \}$, $\bar{V}^* =
\{\bar{\xi}\in V^*|\bar{\xi}=(\xi_1,\cdots,\xi_m,0,\cdots,0)\}$  and
$\tilde{F}^*=F^*|_{\bar{V}^* }$. Then
$\mathcal{L}(\textbf{n})=d\rho\in\bar{V}^*$, where
$\mathcal{L}:(V,F)\to (V^*,F^*)$ is the Legendre transform. Further,
we have
$$ \begin{array}{cccc}
   (V,F) \to&(V^*,F^*)\to  (\bar{V}^*,\tilde F^*)\to &(\bar{V},\tilde F)   \\
   \textbf{n}~~~~\longmapsto & \mathcal{L}(\textbf{n})=d\rho~~~~\longmapsto & \tilde{\mathcal{L}}^{-1}(d\rho)=\tilde{\textbf{n}}
 \end{array}
$$
where $\tilde{\mathcal{L}}:(\bar{V},\tilde F)\to(\bar{V}^*,\tilde
F^*)$ is the Legendre transform. Note that
$\mathcal{L}(\textbf{n})=\tilde{\mathcal{L}}(\tilde{\textbf{n}})=d\rho\in\bar{V}^*$
and $\tilde{F}^*=F^*\big|_{\bar{V}^*}$. Then we have
\begin{equation}\label{5.a}
\left\{\begin{aligned}
      &\tilde F^{*}(d\rho)=F^{*}(d\rho)=\tilde a(\rho)=1,\\
      &\tilde g^{*\alpha\beta}(d\rho)\rho_{\alpha\beta}=g^{*\alpha\beta}(d\rho)\rho_{\alpha\beta}=g^{*ij}(d\rho)\rho_{ij}=\tilde b(\rho).
\end{aligned}\right.\end{equation}
 This means that $\rho$ is also an isoparametric function in $(\bar{V},\tilde F)$. From (\ref{3.4}) we know that $\tilde b(\rho)=-(m-1)k.$ Let $\bar N=\{\bar x\in\bar{V}|\rho(\bar{x})=0\}$ and $\{\tilde{k}_\alpha\}_{\alpha=1}^{m-1}$ be the
constant principal curvatures of $\bar N$. By (\ref{5.a}), we have
$$\sum_{\alpha=1}^{m-1}\tilde{k} _\alpha=\sum_{a=1}^{m-1}k_a=(m-1)k.$$
On the other hand, it follows from  (\ref{3.4}) and (\ref{3.8}) that
$$\sum_{\alpha=1}^{m-1}\tilde{k} _\alpha^2=\sum_{a=1}^{m-1}k_a^2=(m-1)k^2.$$
Therefore, $\tilde{k} _1=\cdots=\tilde{k} _{m-1}=k$. This means that
$\bar N$ is an isoparametric hypersurface with $g=1$ in
$(\bar{V},\tilde F)$. By Theorem 5.3, we conclude that $\bar N$ must
be a (reverse) Minkowski hypersphere $\hat
S_{\tilde{F}_\pm}^{m-1}(\bar x_0,r)$ in $(\bar{V},\tilde F)$ and
thus $N$ must be an (reverse) $F^*$-Minkowski cylinder $\hat
S_{\tilde{F}_\pm}^{m-1}(\bar x_0,r)\times \mathbb{R}^{n-m}$.
\endproof

\section{ \textbf{Isoparametric hypersurfaces in a Randers space}}
Let $(M,F)$ be an $n-$dimensional Randers space, where
$F=\alpha+\beta=\sqrt{a_{ij}y^iy^j}+b_iy^i$.  It is known that
\begin{align}
&g_{ij}=\frac{F}{\alpha}(a_{ij}-\alpha_{y^{i}}\alpha_{y^{j}})+F
_{y^i}F_{y^j}.\label{6.01}\end{align}

 \begin{lem}\label{lem6-1} Let $(M,F)$ be an $n-$dimensional Randers space, where
$F=\alpha+\beta$ and  $b=||\beta||_\alpha$.  Then the Cartan
curvature satisfies $1-Q_{y}(X,Y)=\alpha(1-b^2)>0$ and
$C_y(X,Y,Z)=0$  for any $y,X,Y,Z\in T_xM, x\in M$, where $y,X,Y,Z $
are all orthogonal to each other with respect to $g_y$.
\end{lem}
\proof
  From (\ref{6.01}), we obtain  by a straightforward calculation that
\begin{align} 2C_{ijk}&=\frac{1}{\alpha^2}(\alpha b_k-\beta\alpha_{y^k})h_{ij}
-\frac{F}{\alpha}(\alpha_{y^iy^k}\alpha_{y^j}+\alpha_{y^i}\alpha_{y^jy^k})
+F_{y^iy^k}F_{y^j}+F_{y^i}F_{y^jy^k}, \label{6.0}\end{align}
\begin{align} 2C_{ijkl}=&\left(\frac{F}{\alpha}\right)_{y^ky^l}h_{ij}-\frac{1}{\alpha^2}(\alpha b_k-\beta\alpha_{y^k})(\alpha_{y^iy^l}\alpha_{y^j}+\alpha_{y^i}\alpha_{y^jy^l})\nonumber\\
&-\frac{1}{\alpha^2}(\alpha b_l-\beta\alpha_{y^l})(\alpha_{y^iy^k}\alpha_{y^j}+\alpha_{y^i}\alpha_{y^jy^k})\nonumber\\
&-\frac{F}{\alpha}(\alpha_{y^iy^ky^l}\alpha_{y^j}+\alpha_{y^iy^k}\alpha_{y^jy^l}+\alpha_{y^iy^l}\alpha_{y^jy^k}+\alpha_{y^i}\alpha_{y^jy^ky^l})\nonumber\\
&+F_{y^iy^ky^l}F_{y^j}+F_{y^iy^k}F_{y^jy^l}+F_{y^iy^l}F_{y^jy^k}+F_{y^i}F_{y^jy^ky^l},
\label{6.1}\end{align} where
$h_{ij}=a_{ij}-\alpha_{y^i}\alpha_{y^j}$.
 Let $\{e_{i}\}_{i=1}^{n}$ be an  orthonormal frame field with respect to $g_y$ such that $e_n={\ell}$ and set $e_a=u_a^i\frac{\partial}{\partial x^i}$. Then $$F_{y^i}u_a^i=0,~~~\alpha_{y^i}u_a^i=-\beta(e_a),~~~F_{y^iy^j}u_a^iu_b^j=\alpha_{y^iy^j}u_a^iu_b^j=\delta_{ab},$$
 $$\alpha_{y^iy^jy^k}u_a^iu_b^ju_b^k=F_{y^iy^jy^k}u_a^iu_b^ju_b^k=2C_{ijk}u_a^iu_b^ju_b^k,~~~h (e_a,e_b)=\alpha\delta_{ab},$$
 when $y={\ell}\in S_xM$.
 Thus, we have from (\ref{6.0}) that
\begin{align}
2\hat C_{abc}=2C
(e_a,e_b,e_c)=&\frac{1}{\alpha}\left((\beta(e_c)\alpha
+\beta\beta(e_c))\delta_{ab}
+\delta_{ac}\beta(e_b)+\delta_{bc}\beta(e_a)\right)\nonumber\\
=&\frac{1}{\alpha}\left(\beta(e_c)\delta_{ab}+\beta(e_b)\delta_{ac}+\beta(e_a)\delta_{bc}\right).\nonumber
\end{align}
 It is obvious that
 \begin{align}
\hat C_{abc}&=0,~~~~~~~~~~~~~~~~~~~~a\neq b\neq c\neq a,\nonumber\\
\hat C_{aab}&=\frac{1}{2\alpha}\beta(e_b),~~~~~~~~a\neq b\nonumber\\
\hat C_{aaa}&=\frac{3}{2\alpha}\beta(e_a).\nonumber
\end{align}
By (\ref{6.1}), for any $a\neq b$, we have
\begin{align}
2\hat {\mathcal{C}}_{abab}=\mathcal{C} (e_a,e_b,e_a,e_b)
=&\frac{1}{\alpha^2}(\beta(e_a){\alpha}+\beta\beta(e_a))\beta(e_a)+\frac{1}{\alpha^2}(\beta(e_b){\alpha}+\beta\beta(e_b))\beta(e_b)\nonumber\\
&+\frac{1}{\alpha}\left(2\hat C_{aab}\beta(e_b)-1+2\beta(e_a)\hat C_{bab}\right)+1\nonumber\\
=&\frac{2}{\alpha^2}(\beta(e_a)^2+\beta(e_b)^2)-\frac{1}{\alpha}+1.\nonumber
\end{align}
\begin{align}
4\sum_{c} \hat C_{aac}\hat C_{bbc}=&4( \hat C_{aaa}\hat C_{bba}+\hat C_{aab}\hat C_{bbb})+4\sum_{c\neq a,b}\hat C_{aac}\hat C_{bbc}\nonumber\\
=&\frac{3}{\alpha^2}(\beta(e_a)^2+\beta(e_b)^2)+\sum_{c\neq
a,b}\frac{1}{\alpha^2}\beta(e_c)^2.\nonumber
\end{align}
\begin{align}
1-Q (e_a,e_b)=&1+2\mathcal C_{aabb}-4\sum_{c} \hat C_{aac}\hat C_{bbc}\nonumber\\
=&2+\frac{2}{\alpha^2}(\beta(e_a)^2+\beta(e_b)^2)-\frac{1}{\alpha}-\frac{3}{\alpha^2}(\beta(e_a)^2+\beta(e_b)^2)
-\frac{1}{\alpha^2}\sum_{c\neq a,b}\beta(e_c)^2\nonumber\\
=&2-\frac{1}{\alpha}-\frac{1}{\alpha^2}(|\beta|^2_{g ^*}-\beta(e_n)^2)\nonumber\\
=&1-\frac{\beta}{\alpha}+\frac{\beta}{\alpha^2}-\frac{|\beta|^2_{g
^*}}{\alpha^2}.\nonumber
\end{align}
Using the formula in \cite{BCS}, we have $$|\beta|^2_{g ^*}=g
^{ij}b_ib_j=b^2(\alpha+\beta^2-2\alpha\beta)+\beta^3
$$
Hence
\begin{align}
1-Q (e_a,e_b)=\alpha(1-b^2)>0.\nonumber
\end{align}
\endproof
From Lemma 6.1 and Theorem 5.1 $\thicksim$ 5.5, we have following
classification theorem.

 {\bf Theorem 6.1} In an $n$-dimensional Randers-Minkowski space $(V,F,d\mu)$ where $d\mu$ is
the BH(resp. $HT$)-volume form,  any isoparametric hypersurface has
at most two distinct principal curvatures and it must be either  a
Minkowski hyperplane, a (reverse)  Minkowski hypersphere $\hat
S_{F_{\pm}}^{n-1}(r)$, or an (reverse) $F^*$-Minkowski cylinder
$\hat S_{\tilde{F}_\pm}^{m-1}(r)\times \mathbb{R}^{n-m}$.

Let  $F=\alpha+\beta=\sqrt{\delta_{ij}y^iy^j}+b_iy^i$ be a
Minkowski-Randers metric and $d\mu$ be the BH(resp. $HT$)-volume
form. Set $\lambda=1-b^2,~ b^i=a^{ij}b_j. $ The dual metric $F^*$
can be expressed by(\cite{SZ})
\begin{align}
F^*(\xi)=\alpha^*(\xi)+\beta^*(\xi)=\sqrt{a^{*ij}\xi_i\xi_j}+b^{*i}\xi_i,\quad
\xi=\xi_idx^i \in T^*_x M,\label{4.1} \end{align} where
\begin{align}a^{*ij}=\frac{1}{\lambda^2}(\lambda \delta^{ij}+b^{i}b^{j}),\quad b^{*i}=-\frac{b^{i}}{\lambda}.\label{4.2}\end{align}
Set  $\bar{V}^* := \{\bar{\xi}\in
V^*|\bar{\xi}=(\xi_1,\cdots,\xi_m,0,\cdots,0)\}$  and
$\tilde{F}^*=F^*|_{\bar{V}^* }$.  It is obvious that
$$\tilde{F}^*(\bar{\xi})=\sqrt{\tilde a^{*{\alpha }{\beta }}\xi_{\alpha }\xi_{\beta }}+\tilde b^{*{\alpha }}\xi_{\alpha }=\frac{1}{\lambda}\left(\sqrt{(\lambda \delta^{{\alpha }{\beta }}
+b^{{\alpha }}b^{{\beta }})\xi_{\alpha }\xi_{\beta }} -b^{{\alpha
}}\xi_{\alpha }\right).$$ Therefore $$\tilde a^*_{{\alpha }{\beta
}}=\lambda\left(\delta_{{\alpha }{\beta}}-\frac{1}{\lambda+\bar
b^2}b_{\alpha }b_{\beta}\right),~~~~ \tilde b^*_{{\alpha }}=\tilde
a^*_{{\alpha }{\beta }}\tilde{b}^{*{\beta }}=-\frac{\lambda
}{\lambda+\bar b^2}b_{\alpha }, $$ where $\bar b^2=b^{\alpha
}b_{\alpha }$. The dual metric of $\tilde{F}^*$ in $\bar V=\{y\in
V|y=(y^1,\ldots,y^m,0,\ldots,0) \}$ can be written as
$$\tilde{F}(\bar{y})=\frac{1}{\bar\lambda}\left(\sqrt{(\bar\lambda \tilde a^*_{{\alpha }{\beta }}
+\tilde b^*_{{\alpha }}\tilde b^*_{{\beta }})y^{\alpha
}y^{\beta}}-\tilde b^*_{{\alpha }}y^{\alpha }\right)=\sqrt{\tilde
a_{{\alpha }{\beta}} y^{\alpha }y^{\beta }}+\tilde b_{{\alpha
}}y^{\alpha },$$ where $\bar\lambda=1-\tilde{b}^{*\alpha
}\tilde{b}^*_{\alpha}=1-\frac{\bar b^2}{\lambda+\bar b^2}$ and
$$\tilde a_{{\alpha }{\beta }}=\frac{1}{\bar\lambda^2}\left(\bar\lambda \tilde a^*_{{\alpha }{\beta }}
+\tilde b^*_{{\alpha }}\tilde b^*_{{\beta }}\right)=(\lambda+\bar
b^2)\delta_{\alpha\beta},~~\tilde b_{\alpha
}=-\frac{1}{\bar\lambda}\tilde b^*_{\alpha }= b_{ \alpha}.$$ That is
$$\tilde{F}(\bar{y})=\sqrt{(\lambda +\bar b^2)}|\bar y|+ b_{{\alpha
}}y^{\alpha }.$$

{\bf Example 3} Let
$F=\alpha+\beta=\sqrt{\delta_{ij}y^iy^j}+b_iy^i=|y|+\beta(y)$. The
equation of the isoparametric cylinder $\hat
S_{\tilde{F}_\pm}^{m-1}(r)\times \mathbb{R}^{n-m}$  can be expressed
explicitly as
$$
\sqrt{\lambda+\bar b^2}|\bar x|\pm\beta(\bar x)=r,~~~~ \bar
x=(x^1,\ldots,x^m,0,\ldots,0).$$

 From \cite{SZ}, the gradient vector of $f$ at $x$ can be written as
 \begin{align}
\nabla f=\frac{F^*(df)}{\lambda
\alpha^*(df)}(f^i-F^*(df)b^i)\frac{\partial}{\partial x^i},
\label{4.4}\end{align} where $df=f_idx^i,f^i=\delta^{ij}f_j$. Let
$\nabla^\alpha f$ and $\Delta_\sigma^\alpha f$ be the gradient
vector and the Laplacian of $f$ with respect to $d\mu$,
respectively. Then
 \begin{align}
\nabla f=&\frac{F^*(df)}{\lambda \alpha^*(df)}(\nabla^\alpha f-F^*(df)\beta^{\sharp}),\label{4.5}\\
\Delta f=&\text{div}\nabla f=\frac{F^*(df)}{\lambda
\alpha^*(df)}\Delta^\alpha f+<df-F^*(df)\beta,
d\left(\frac{F^*(df)}{\lambda \alpha^*(df)}\right)>_\alpha\nonumber\\
&-\frac{F^*(df)}{\lambda \alpha^*(df)}<dF^*(df),
\beta>_\alpha,\label{4.6}\\
F^*(df)=&\alpha^*(df)+\beta^*(df)=\frac{1}{\lambda}\left(\sqrt{\lambda|df|_\alpha^2
+<df,\beta>_\alpha^2}-<df,\beta
>_\alpha\right).
\label{4.7}\end{align} (\ref{3.1}) can be written as
 \begin{equation}\label{4.8}\left\{\begin{aligned}
      &|df|_\alpha^2=\lambda \tilde a ^2+2\tilde a \zeta,\\
      &\Delta^\alpha f=\frac{\lambda \alpha^*\tilde b }{\tilde a }-\lambda \tilde a  \tilde a '+\frac{1}{\lambda \alpha^*}<df-\tilde a \beta,
      d(\lambda \alpha^*)>_\alpha,
\end{aligned}\right.\end{equation}
where $\zeta=<df,\beta>_\alpha.$ Since $\lambda \alpha^*=\lambda
\tilde a +\zeta $, we have
\begin{align}d(\lambda \alpha^*)=&\lambda \tilde a'df+d\zeta,\nonumber\\
 <df,d\zeta>_\alpha=&<<\nabla^\alpha df,\beta>_\alpha,df>_\alpha
=<<\nabla^\alpha df,df>_\alpha,\beta>_\alpha. \label{4.9}\end{align}
 It follows from $(\ref{4.8})_1$ that
\begin{align*}
 <\nabla^\alpha df,df>_\alpha=&(\lambda \tilde a +\zeta)\tilde a 'df+\tilde a  d\zeta,\nonumber\\
 <<\nabla^\alpha df,df>_\alpha,\beta>_\alpha=&(\lambda \tilde a +\zeta)\tilde a '\zeta+ \tilde a <\beta,d\zeta>_\alpha.
\end{align*}
Thus
\begin{align*}<df-\tilde a \beta,d(\lambda \alpha^*)>_\alpha=<df-\tilde a \beta,d\zeta>_\alpha+\lambda\tilde{a}'(|df|^2_\alpha-\tilde{a}\zeta)=\tilde a '(\lambda \tilde a +\zeta)^2.
\end{align*}
 \begin{thm}\label{thm6-4}Let $(V,F,d\mu)$ be an $n$-dimensional Randers-Minkowski space where $d\mu$ is
the BH(resp. $HT$)-volume form.  A function $f$ is isoparametric if
and only if it satisfies
\begin{equation}\label{4.11}\left\{\begin{aligned}
      &|df|_\alpha^2=\lambda \tilde a (f)^2+2\tilde a (f)<df,\beta>_\alpha,\\
      &\Delta^\alpha f=\lambda \tilde b (f)+\left(\frac{\tilde b (f)}{\tilde a (f)}+\tilde a '(f)\right)<df,\beta>_\alpha.
\end{aligned}\right.\end{equation}
 \end{thm}

{\bf Example 4} Let $f(x)=|\bar x|+\beta(x)$, where $\bar
x=(x^1,\ldots,x^m,0,\ldots,0)$. Then
$$f_{\alpha }=\frac{x^{\alpha }}{|\bar x|}+b_{\alpha },~~~~f_{{\lambda }}=b_{{\lambda }},~~~~f_{{\alpha }{\alpha }}=\frac{|\bar x|^2-x^{\alpha } x^{\alpha }}{|\bar x|^3},~~~~f_{\lambda {\lambda }}=0,
$$
 $$ |df|_\alpha^2=1+{b}^2+2\frac{\beta(\bar x)}{|\bar x|},~~~~<df,\beta>_\alpha=\frac{\beta(\bar x)}{|\bar x|}+{b}^2,~~~\Delta^\alpha f=\frac{m-1}{|\bar x|}.$$
  By choosing $\tilde a (t)=1$, one can easily to check that
 $$\lambda \tilde a (f)^2+2\tilde a (f)<df,\beta>_\alpha=\lambda+2\left(\frac{\beta(\bar x)}{|\bar x|}+b^2\right)= |df|_\alpha^2.$$
 That is, $f$ is a transnormal function. If $f$ is an isoparametric function, then from (\ref{3.8}), its each regular level surface in $M_{f}$  must have two constant principal curvatures $0$  and $k=\frac{-1}{s+c}$. Denote the  multiplicity of $k$ by $q$. It follows from (\ref{3.4}) that $\tilde b (t)=\frac{q}{t+c}$, where $c$ is a constant. Then when  $b_{{\lambda }}\neq 0$,
 \begin{align*}\lambda \tilde b (f)+\left(\frac{\tilde b (f)}{\tilde a (f)}+\tilde a '(f)\right)<df,\beta>_\alpha&=\frac{q}{|\bar x|+\beta( x)+c}\left(\lambda+\frac{\beta(\bar x)}{|\bar x|}+b^2\right)\\
 &=\frac{q(|\bar x|+\beta(\bar x))}{|\bar x|(|\bar x|+\beta( x)+c)}\neq\Delta^\alpha f
 \end{align*}
 at points where $\beta(x)\neq\beta(\bar x)$.
That is, $f$ does not  satisfy (\ref{4.11}) and thus it is not an
isoparametric function.

\bibliographystyle{Plain}

\end{document}